\documentclass[a4paper]{amsart}
\usepackage[latin1]{inputenc}
\usepackage{amsmath}
\usepackage{amsthm}
\usepackage{latexsym}
\usepackage{amssymb}
\usepackage{wasysym}
\usepackage{setspace}
\usepackage{mathrsfs}
\usepackage[]{todonotes}
\usepackage{slashbox}
\usepackage{paralist}
\usepackage{booktabs}
\usepackage{tabularx}

\makeatletter

\def\R{\ifmmode{\mathbb R }\else{$\mathbb R $ }\fi}
\def\Q{\ifmmode{\mathbb Q }\else{$\mathbb Q $ }\fi}
\def\N{\ifmmode{\mathbb N }\else{$\mathbb N $ }\fi}
\def\Rn{\ifmmode{\mathbb R^n }\else{$\mathbb R^n $ }\fi}
\def\C{\ifmmode{\mathbb C }\else{$\mathbb C $  }\fi}
\def\F{\ifmmode{\mathbb F }\else{$\mathbb F $ }\fi}
\def\D{\ifmmode{\mathscr{D} }\else{$\mathscr{D} $ }\fi}
\newcommand{\norm}[1]{\left\Vert #1 \right\Vert}

\renewcommand{\limsup}{\overline{\lim}}
\renewcommand{\liminf}{\underline{\lim}}
\newcommand{\Lloc}{L_{\text{loc}}^1}
\renewcommand{\phi}{\varphi}
\newcommand{\Wloc}{W_{\text{loc}}^{1,n}}

\newcommand{\abstand}{\vspace{0.8cm}}

\newtheorem{definition}{Definition}[section]
\newtheorem{theorem}[definition]{Theorem}
\newtheorem{proposition}[definition]{Proposition}
\newtheorem{corollary}[definition]{Corollary}
\newtheorem{example}[definition]{Example}
\newtheorem{examples}[definition]{Examples}
\newtheorem{lemma}[definition]{Lemma}
\newtheorem{remark}[definition]{Remark}

\allowdisplaybreaks
\numberwithin{equation}{section}

\begin{document}
\title{Spaces of rapidly oscillating functions}
\author{Hans Günzler}
\begin{abstract}
For various function spaces of the form $gU$ or $U + gV$, $U$ and $V$ e.g. AP, BUC or UC, $g(t):=e^{i |t|^r}$, their properties are discussed, especially a Loomis type condition ($\Delta$) which has been essential in getting such properties: $A$ satisfies ($\Delta$) if $f \in L_{\text{loc}}^1(J, X)$ and $f (\ \cdot + h) - f( \cdot) \in A$ for all $h>0$ implies $f - \frac{1}{h} \int_0^h f(\ \cdot +s ) ds \in A$.
\end{abstract}
\keywords{Loomis conditions, rapidly oscillating, oscillatory conditions, Porada inequality}
\subjclass[2010]{46E40, 46E99, 42A75}
\maketitle
\tableofcontents
\newpage
%

\section{Introduction}
The spaces $A\subset X^J$, $J=$ halfline or \R, $X$ a Banach Space, of ``rapidly oscillating functions'' we are interested in are of the form $\gamma U$ or $U+\gamma V$, with e.g. $\gamma(t)=e^{i \omega |t|^r}, 0\ne \omega \in \R, 1<r$ and $U,V$ general function spaces of e.g. almost periodic functions AP or (bounded) uniformly continuous functions (B)UC.\\
We discuss properties of such $A$ as translation invariance, closedness with respect to uniform convergence on $J$, when one has $A\subset$ mean extension $MA$, meaning mollifier $M_h(f):= \frac{1}{h}\int_0^h f(\cdot + s) ds \in A$ if $f \in A, h>0$, and especially a condition ($\Delta$), which has been introduced in \cite[Def.1.4 p.7]{lit5}, \cite[p. 677]{lit7}, meaning if $f\in L_{\text{loc}}^1(J,X) $ with all differences $\Delta_hf\in A, (\Delta_hf)(t):=f(t+h)-f(t)$, then $f-M_h(f) \in A, h>0$.\\
This condition has been found useful in many situations, e.g.
\begin{itemize}
\item for getting Loomis' condition ($L_b$) (Definition (6.4) below, needed in the study of asymptotic properties of solutions of differential-difference systems, see e.g. \cite[Lemma 2.4]{lit4} from the Bohl-Bohr-Kadets Theorem (P$_b$) ( (6.5), \cite[Prop.3.12  p. 682]{lit7})
\item for an explicit description of the class of functions $f$ whose differences $\Delta_hf$ are in a given class $A$ ( (2.9), Prop. 4.7, (4.19))
\item similarly for an explicit description of the mean extension MA (Prop.3.8)
\item showing that these MA are closed under multiplication with $e^{i\omega t}, \omega \in \R$, if $A$ is (\cite[Prop.8.1 p.52]{lit6})
\item generalizing Esclangon-Landau results to solutions of differential-difference systems (\cite[Theorem 5.1 and Corollaries]{lit8})
\item in the proof of our main result, ($\Delta$) for BUC+ $\gamma$BUC, the condition ($\Delta$) is needed 3 times: for BUC, for the $W$ of (6.29) and for BUC($\R$,X)
\end{itemize}
Main results of the present note are:
\begin{itemize}
\item  $\gamma (B)UC\subset M(\gamma (B)UC)$: Prop. 3.11, Examples 3.16, 3.17
\item  a Porada type inequality for $U+\gamma V$ (Prop. 5.4) and completeness of such $U+\gamma V$ (Theorem 5.6)
\item ($\Delta$) for $\gamma BUC$, $\gamma UC$ (Theorem 4.3), $BUC+\gamma BUC$ (Theorem 6.5, Cor. 6.6) and related vector sums ( ($\Delta$)-table p. 25, table p. \pageref{tabelle})
\item ($P_b$),($L_b$) for $(B)UC + \gamma (B)UC$ (Corollaries 6.9, 6.13)
\end{itemize}
What is needed for all this is collected in §§2, 3, 4; these results are mostly due to Bolis Basit and the author (\cite{lit3}-\cite{lit12}), cited here, occasionally in expanded form.\\
The autor wants to thank Bolis Basit for this fruitful collaboration, for his ideas, perseverance and hard work.
\abstand
\section{Notation, Definitions, preliminary Lemmas}
In the following \R  resp. \C denotes the real resp. complex field, $\R_+:=[0,\infty), \R^+:= (0,\infty), J$ will always be an interval $\subset \R$ of the form $[\alpha, \infty)$, $(\alpha, \infty), \R, \alpha \in \R ; \N:=\{1,2,\dots\}, \N_0:=\{0\} \cup \N$.\\
$X$ will denote a Banach space, nonempty and $\ne \{0\}$, with scalar field $\F=\F_X \in \{\R,\C\}$ and norm $\norm{ \cdot }$.\\
For $f:J\rightarrow X $ the translate $f_a$, the difference $\Delta_af$ and $|f|:J\rightarrow \R_+$ are defined by, $t\in J, a\in \R_+$ resp. \R if $J=\R$,
\begin{equation}
f_a(t)=f(t+a), \quad \Delta_af=f_a-f, \quad |f|(t)=\norm{f(t)}.
\end{equation}
For $\omega \in \R, 0\ne r\in \R, t\in \R$
\begin{equation}
\gamma_{\omega}(t):=e^{i\omega t}
\end{equation}
\begin{equation}
g_{\omega, r}(t):=e^{i\omega |t|^r},g(t):=g_{1,2}(t)=e^{it^2};
\end{equation}
wherever the $g_{\omega,r}$ or $\gamma_{\omega}$ appear, $\F_X=\C$ is assumed.
\begin{equation}
\limsup_{|t|\to \infty}:=\limsup_{t\to\infty} \text{ if }J\ne\R, \limsup_{|t|\to \infty}:=\max \{\limsup_{t\to \infty}, \limsup_{t\to-\infty}\} \text{ if } J=\R
\end{equation}
similarly for $\liminf_{|t|\to\infty}$, with ``min''.\\
All function spaces are $\subset X^J$, with pointwise defined $=,+,\beta \ \ (\beta \in \F)$.\\
For $k\in \N_0, \  C^k(J,X)$ contains all $k$ times continuously differentiable functions $f:J\rightarrow X, \ C(J,X):=C^0(J,X)$.\\
$X_c, C_0(J,X), BC(J,X),UC(J,X),BUC(J,X), AP(J,X)$ contain all continuous\\ $f:J\rightarrow X$ which are\\
constant, vanishing at infinity and uniformly continuous, bounded, uniformly continuous, bounded uniformly continuous resp. almost periodic (=ap), for ap on \R see e.g. \cite[p.3]{lit1},\cite[p.1]{lit24},\cite[Theorem 4.1]{lit12}, $AP(J,X):=AP(\R,X)\big|J$.\\
The class of asymptotic ap functions is defined by $AAP(J,X):=C_0(J,X)+AP(J,X)$ [29, p. 35 p. 46]. For further generalizations $BAA,AA,LAP,REC, PS_+$ see (3.6) and the reference after (3.6); for spaces of ergodic functions $E, BE, EM,$\\ $ EM_0, BEM,BEM_{(0)}, CEM_{(0)},$ see (4.22) and after, for $Av, Av_0$ see (6.2), (6.3), for $PAP, AAA,$ $PAA, EAP$ after (6.3).\\
If $f\in X^J, \norm{f}$ always means $\norm{f}_{\infty}:=\sup\{\norm{f(t)}: t\in J\}$, for the Stepanoff norm $\norm{ \cdot }_{S^1}$ see (3.4).\\
$L_{\text{loc}}^1(J,X)$ contains all $f:J\rightarrow X$ which are Bochner-Lebesgue integrable over every compact intervall $\subset J, \  L^p(J,X)$ contains all Bochner-Lebesgue measurable $f:J\rightarrow X$ for which $|f|^p$ is Lebesgue-integrable over $J$ if $1\leq p<\infty$, resp. $f$ is bounded everywhere on $J$ if $p=\infty$, with corresponding seminorms (all linear subspaces of $X^J$, no equivalence classes). The Sobolev spaces $W_{\text{loc}}^{1,n}$ are defined in Example 3.6.\\
For $U,V \subset X^J$, $U+V:=$ vector sum $\{u+v: u\in U, v\in V\}.$ For $f\in L_{\text{loc}}^1(J,X)$, with fixed $\alpha_0 \in J, \alpha_0=0 $ if $0\in J$, we use
\begin{equation}
Pf: J\rightarrow X, (Pf)(t):=\int_0^t  f(s)ds, \ t \in J
\end{equation}
(all integrals in the following are Bochner-Lebesgue integrals, usually for $X-$valued functions).\\
For any $A\subset X^J$
\begin{equation}
A \textit{ positive invariant means } f_a\in A \text{ for all } f\in A, a\in \R_+,
\end{equation}
\begin{equation}
A \textit{ invariant means } f_a\in A \text{ for all } f\in A, a\in \R \text{ with } J=\R,
\end{equation}
\begin{align}
A \textit{ uniformly closed: if  } f_n\in A, n\in \N, f\in X^J, f_n (t) \to f(t) \text{ as } n\to \infty,\\ \text{ uniformly in } t\in J, \text{ implies } f\in A. \nonumber
\end{align}
For $f\in \Lloc(J,X)$ and $h\in \R^+$, the \textit{means} (mollifiers) $M_hf $ are defined by
\begin{equation}
(M_h f)(t):=\frac{1}{h}\int_0^h f(s+t)ds=\frac{1}{h}\int_t^{t+h} f(s) ds, \ t\in J;
\end{equation}
one has
\begin{equation}
M_h(f_a)=(M_hf)_a, a\in \R_+ \text{ resp. } a\subset \R \text{ if } J=\R
\end{equation}
The \textit{mean extensions} $M^nA$, introduced in \cite[p. 120]{lit4},\cite[Def.1.1]{lit5} are defined by
\begin{align}
MA:=M^1A:=\{f\in\Lloc (J,X):M_hf \in A \text{ for all } h \in \R^+\},\\
	M^{n+1}A:=M(M^nA), M^0A:=A\cap \Lloc, \  n \in  \N,\nonumber
\end{align}
one has $M(M^0A)=MA$; \\
usually $A\subset MA$, then $A\subset MA \subset M^2A \subset \dots$ (see §3 (3.1), Example 3.4).\\

\textit{Difference classes} $\Delta^n A$ (introduced in \cite[p. 680]{lit7} ) are, for $A \subset X^J$, defined by
\begin{align}
\Delta A:= \Delta^1 A:=\{f\in \Lloc(J,X): \Delta_h f \in A \text{ for all } h\in \R_+\},\\
\Delta^{n+1} A:= \Delta(\Delta^n A), \  n \in \N. \nonumber
\end{align}
An $A\subset X^J$ \textit{satisfies} ($\Delta$) (introduced in \cite[p. 677]{lit7}, see the introduction) means
\begin{align}
\text{if } f \in \Lloc (J,X) \text{ with all } \Delta_af \in A, a\in \R^+, \text{ then } f-M_hf \in A \\ \text{ for all }h\in R^+;\nonumber
\end{align}
\begin{align}
A \textit{ satisfies } (\Delta_1) \text{ means (2.13) with ``for all } h\subset \R^+\text{''} \\ \text{replaced by ``for at least one } h \in R^+ \text{''}\nonumber
\end{align}
The oscillatory conditions $O_1$ and $O_2$ are defined in §5, (5.3) and (5.10). ($\Delta P$) and $(\Gamma)$ are defined in § 8.
\begin{lemma}
$I$ arbitrary interval $\subset \R, \ \epsilon_0>0, n\in \N_0, f\in \Lloc (I,X);$
if then $\Delta_hf \in C^n(I^{-h},X)$ for all $h\in (0,\epsilon_0)$, then $f \in C^n(I,X)$.

Here $I^{-h}:=(\alpha, \beta -h]$ resp. $(\alpha,\beta-h]$ if $I=(\alpha,\beta)$ resp. $(\alpha,\beta],$  similarly for $[\alpha,\beta)$ resp. $[\alpha,\beta]$, $ \ -\infty\leq \alpha < \beta \leq \infty$.
\begin{proof}
\cite[Prop 1.3 p. 1009/1010]{lit8}, \cite[Prop 1.5 p. 778]{lit5}.
-- \end{proof}
\end{lemma}

\begin{lemma}
If $I, \epsilon_0, n,f$ are as in Lemma 2.1, line 1, and if to each $h \in (0, \epsilon_0)$ there exists $\phi \in C^n(I^{-h},X)$ with $\Delta_hf=\phi$ a.e. on $I^{-h}$, then there exists $\Phi \in C^n(I,X)$ with $f=\Phi$ a.e. on $I$.
\begin{proof}
This follows from \cite[Lemma 4.4 p. 30]{lit6}, a distribution version of Lemma 2.1 above. A direkt proof using only Lemma 2.1 for $M_hf $ is possible, we omit this.
-- \end{proof}\end{lemma}
For generalizations of Lemma 2.1, without $f\in \Lloc$, see § 7.
\abstand
\section{Mean classes}
By definition (2.11), always $MA \subset \Lloc(J,X)$; if $U\subset V\subset X^J, MU \subset MV$; thus $A \subset MA$ implies $M^nA\subset M^{n+1}A, n\in\N; M_h(f_a)=(M_hf)_a.$ In (2.11) under suitable assumptions fewer than all $h\in \R^+$ suffice, for this see \cite[p. 24 §3 (a), (b), (c)]{lit6}.
\begin{proposition}
If $A\subset \Lloc (J,X)$ is linear, positive invariant and satisfies $(\Delta)$, then $A\subset MA$.
\begin{proof}
By the assumptions, if $f\in A$, then $\Delta_hf\in A$, with $(\Delta)$ then $f-M_hf \in A$, so $M_hf\in A, h\in \R^+$.
-- \end{proof}
\end{proposition}

By example 3.14, $A=gAP$, Proposition 3.1 becomes false without ($\Delta$), even if $A$ is additionally uniformly closed and $\subset BUC(\R,X)$; by example 4.4, Proposition 3.1. becomes false without $A$ positive invariant.

\begin{proposition}
If $A$ is convex (e.g. linear), positive invariant, uniformly closed and $\subset  UC(J,X), $ then
\begin{equation}
A \subset MA \subset M^2 A \subset \dots
\end{equation}
\normalfont For a generalization see \cite[Prop. 3.2 p . 25/26]{lit6}
\begin{proof}
\cite[Prop. 2.2 p. 1011]{lit8}. -- \end{proof}
\end{proposition}

\begin{proposition}
If $n \in \N$ and $A\subset \Lloc (J, X)$ is linear resp. positive invariant resp. uniformly closed resp. invariant resp. satisfies ($\Delta_1$) resp. satisfies ($\Delta$) resp. satisfies $A\subset MA$, then $M^nA$ has the same property.
\begin{proof}
\cite[Proposition 2.1 p. 1010]{lit8}. -- \end{proof}
\end{proposition}

\begin{examples}
$C^k(J,X), \ C_0(J,X), \ AP(J,X), \ BUC(J,X), \ UC(J,X), BC(J,X),$ \\ $ L^p(J,X), 1\leq p \leq \infty$, all satisfy (3.1);
for further examples see \cite[Examples 3.4, 3.5, 3.6, 3.13]{lit6}, \cite[Examples 2.3, 2.5]{lit8}.
\begin{proof}
This follows from the properties of $M_h$ resp. Proposition 3.2 resp. the continuous Minkowski inequality \cite[p. 251 Aufgabe 7.92, above (4.2)]{lit20} for $L^p$.
-- \end{proof}
\end{examples}
For examples of $A$ where all the ``$\subset$'' in (3.1) are strict see \cite[Prop. 3.8, Examples 3.9, Remark 3.11]{lit6}

\begin{example}
$M^nX_c=X_c  + \{f\in X^J:f=0 \text{ a.e.}\}$ (Fubini or Lemma 2.2)
\end{example}

\begin{example}
If $A=$ Sobolev space $W_{\text{loc}}^{1,n}(J,X):=\{f\in C^{n-1}(J,X): f^{(n-1)}$ locally absolutely continuous on $J$ and $(f^{(n-1)})'$ exists a.e. in $J\}, n\in \N$, then
\begin{align}
A\subset MA \subset \dots \subset M^nA=M^{n+1}A=\dots, M\Lloc (J,X)=\Lloc (J, X), \\
\text{with all ``$\subset$'' strict.}\nonumber
\end{align}
(For $\Wloc$ see \cite[Propositions 1.2.2 and 1.2.3]{lit2}, \cite[Theorem 3.8.6 p. 88]{lit22},\cite[(2.5) p. 16]{lit6}.)
\end{example}

\begin{proposition}
(i) If $A$ is linear, positive invariant, uniformly closed $\subset L^{\infty} (\R, X)$ with $A\star \mathscr D(\R,\F) \subset A \subset MA$, then
$$ \mathscr{D}'_A(\R,X) \cap \Lloc(\R,X) \subset \bigcup_0^{\infty} M^nA.$$
(ii) If in addition to the assumptions of (i) $A$ satisfies ($\Delta_1$), then
$$ \mathscr{D}'_A \cap \Lloc =\bigcup_0^{\infty} M^nA.$$
Here $\mathscr{D}(\R,\F):=\{\phi \in C^{\infty}(\R,\F): \text{supp }  \phi \text{ compact}\}$, Schwartz's test functions;
$\mathscr{D}'_A(\R,X):=\{T\in \mathscr{D}'(\R,X): T\star \phi \in A \text{ for all } \phi \in \mathscr{D}(\R,\F)$.
See \cite[Def. 1.2]{lit12}.
\begin{proof}
``$\subset$'': As in \cite[proof of Prop. 3.7]{lit4}, with \cite[Theorem 2.4]{lit12}.\\
``='': \cite[(2.19)]{lit12}.
-- \end{proof}
\end{proposition}
For relations between $A\subset MA$ and $A\star \mathscr{D}\subset A$ see \cite[(2.20)]{lit6}.
Examples for Proposition 3.7 (ii): \cite[Examples 3.2/3/4/5/6]{lit12}.

\begin{proposition}
If $A\subset \Lloc (J,X)$ is linear, positive invariant and satisfies ($\Delta$), then
\begin{equation}
MA=A+A'
\end{equation}
with $A':=\{\phi \in \Lloc(J,X): \text{ to } \phi \text{ exists } \Phi \in A\cap W_{\text{loc}}^{1,1} (J,X) \text{ with } \phi =\Phi' \text{ a.e. on } J\}$
\begin{proof}
\cite[Prop.5.1 (ii) p. 38]{lit6}; for $W_{\text{loc}}^{1,1}$ see Ex. 3.6; for extensions \cite[Prop. 5.1]{lit6}\mbox{.
--} \end{proof}
\end{proposition}

\begin{proposition}
If $ A\subset \Lloc(J,X)$ is uniformly closed, $n\in \N$, then $UC(J,X)\cap M^nA \subset A$.
\begin{proof}
\cite[Prop.2.9 p. 1012]{lit8}.
-- \end{proof}
\end{proposition}
\newpage
In the following we need
\begin{lemma}
If $\gamma \in \Lloc (J,X), a\in J$ and $\lim_{T\to \infty} \int_a^T\gamma(s) ds $ exists $\in X, J\ne \R$, then $\gamma \in MC_0(J,X)$; the same holds for $J=\R$, if additionallay $\lim_{T\to - \infty}\int_T^a \gamma (s)ds$ exists $\in X$;\\
special case: $\gamma=|t|^sg_{\omega,r}$ with $0\ne \omega \in \R, \  0\leq s<r-1\in \R$.
\begin{proof}
The general case follows from the definitions; the case $|t|^sg_{\omega,r}$ with integration by parts ($\int_0^{\infty} g_{1,2}ds=\frac{\sqrt{\pi}}{2} e^{i\pi/4}$ \cite[tables 19.6.1: 20]{lit15}).
-- \end{proof}
\end{lemma}

\begin{lemma}
If $\gamma \in MC_0(J,\F)$ and Stepanoff-norm $\norm{\gamma}_{S^1} < \infty, $ then $\gamma BUC(J,X) \subset MC_0(J,X)$.\\

\textnormal{Here for} $f\in \Lloc (J,X)$.
\begin{equation}
\norm{f}_{S^1}:=\sup \left\{\int_t^{t+1}\norm{f(s)} ds: t\in J\right\}
\end{equation}
(\cite[p. 71]{lit14},\cite[p. 132]{lit4}).
\textnormal{Already} $gUC$ \textnormal{not} $\subset MC_0$, \textnormal{only} $\subset M(gUC)$ \textnormal{(see Examples 3.18, 3.19):}
$gt \notin MC_0$ \textnormal{(but see Prop. 3.16)}.\\
\textnormal{Without} $\norm{\gamma}_{S^1} <\infty$ \textnormal{, at least }$|t|^sg_{\omega,r} AP\subset MC_0$ if $0\leq s<r-1,\ 0 \ne \omega\in \R$ (\textnormal{but} $1\in |t|^sg_{\omega,r}C_0(J,X)$ \textnormal{if} $s>0$).
\begin{proof}
To $f\in UC(J,X)$ and $\epsilon>0$ exists $\delta>0$ with $\norm{f(t+s)-f(t)} \leq \epsilon$ if $0\leq s\leq \delta, t\in J$; then if $h>0, \frac{1}{\delta} \leq n\in \N, h_j:=(jh)/n, $ one has
\begin{align}
 \norm{hM_h(\gamma f)(t)-\sum_{j=1}^n f(t+h_j)\int_{h_{j-1}}^{h_j} \gamma (t+s)ds} \leq \\
\epsilon \int_0^h|\gamma(t+s)|ds \leq \epsilon (h+1) \norm{\gamma}_{S^1}, \text{ all } t\in J. \nonumber
\end{align}
$\gamma \in MC_0$ implies $\phi_j:=\int_{h_{j-1}}^{h_j} \gamma (\cdot \ +s)ds \in C_0, 1\leq j\leq n; $ if $f\in BUC$,
$\Phi_{\epsilon}:=\sum_{1}^n f_{h_j}\phi_j \in C_0$; since $C_0$ is uniformly closed, $M_h(\gamma f)\in C_0$.
-- \end{proof}
\end{lemma}

\begin{proposition}
If $\gamma \in C(J,\F)$ is as in Lemma 3.11 with $\inf_J |\gamma| >0,
C_0(J,X) \subset A\subset BUC(J,X), $ then $ \gamma A \subset M(\gamma A)$.
\begin{proof}
Lemma 3.11 and $C_0=\gamma (\frac{1}{\gamma}C_0)\subset \gamma C_0 \subset \gamma A.$
-- \end{proof}
\end{proposition}
Special such $A, \gamma$: Asymptotic almost periodic functions $AAP:=C_0 + AP$, Eberlein weakly almost periodic functions $EAP,BUC; \ \gamma_{\omega,r}$ with $1<r\in \R, 0\ne\omega \in  \R$ (\cite[section 2]{lit13}, references in \cite[p. 425]{lit9}, \cite[Def. 3.1, Theorem 3.3]{lit29}).
Without $\inf_J|\gamma| >0$ Proposition 3.12 becomes false : $\gamma=\cos (t^2) \in MC_0$ by Lemma 3.10, but already $M_h\gamma$ not in $\gamma BUC$ for small $h>0$ ($t^2=\pi/2).$ See also Example 3.14.\\
Also, $\gamma \in MC_0$ cannot be weakened to $\gamma$ satisfies $O_1$ (see (5.3),(5.4)): $\gamma BUC \not\subset M(\gamma BUC)$ if $\gamma= 2+ g_{1,2}$.\\

Concerning the title of this note, examples of rapidly oscillating functions are the $g_{\omega,r}$ with $1<r\in \R, 0\ne \omega \in \R$, especially $g=g_{1,2}=e^{it^2}$.\\
To get a linear invariant uniformly closed space containing $g$, one has to consider $A_{ro}=\norm{ \quad  }_{\infty}$-closure of linear hull of $\{$translates $g_a: a\in\R\}$; this is precisely $gAP(\R,\C)$ (\cite[Ex. 3.1, p. 25]{lit6}). Unfortunately, this $A$ is in some sense pathological by Example 3.14 below. For this we need first
\begin{lemma}
(a) If $\gamma$ is as in Lemma 3.11 and $ 0 \in A\subset PS_+(J,X)$ of (3.6), then $$A\cap \gamma BUC(J,X)=\{0\}$$
(b) If $\gamma$ is as in (a), $A_1,A_2 \in PS_+(J,X)$ and $\liminf_{|t| \to \infty}|\gamma| >0$, then $$(A_1 + \gamma(A_2\cap BUC)) \cap C_0(J,X)=\{0\}$$
\end{lemma}
Here the class of \textit{Poisson stable} functions \cite[p. 80 Def. 1]{lit24} is defined by
\begin{align}
&PS_+(J,X):=\{f\in C (J,X): \text{ to each } n\in\N \text{ exists } \tau \in [n,\infty) \text{ with }\\ &\norm{f(t)-f(t+\tau)}\leq \frac{1}{n}, \ |t|\leq n, \  t\in J\}. \nonumber
\end{align}
One has $X_c \subset AP \subset \text{ Bochner almost automorphic functions } BAA \subset$\\
$\text{Veech almost automorphic functions } AA \subset \text{ Levitan almost periodic functions } LAP \subset \text{recurrent functions } REC \subset PS_+$
(see \cite[p. 430 (3.3)]{lit9}, $BAA\cap UC=AA\cap UC$).\\
Any  of these $A$'s can be used in Lemma 3.13 and Example 3.14.

Since $\frac{1}{g_{\omega, r} }=g_{-\omega,r}$, with Lemma 3.10 and Lemma 3.13(a) one gets
\begin{equation}
(g_{\omega,r } PS_+(J,X)) \cap BUC(J,X)=\{0\} \text{ if } 0\ne \omega\in \R, 1<r \in \R.
\end{equation}
For a generalization see Corollary 5.2, and also Lemma 5.1.\\
(b) of Lemma 3.13 becomes false without $\liminf |\gamma| >0: \gamma=1/(1+|t|), f_2=1, \gamma \in \gamma PS_+\cap BUC.$\\
(b) becomes trivially false if ``$C_0$'' is replaced by ``$BUC$'': $AP+\gamma \{0\}\subset BUC$.

\textbf{Proof} of (a) of Lemma 3.13: If $f \in PS_+ \cap \gamma BUC, M_hf\in C_0$ by Lemma 3.11; since $PS_+ \subset MPS_+ \ (|M_h(\Phi)|\leq M_h(|\Phi|))$ and $PS_+\cap C_0=0$ by definition (3.6), $M_hf=0$ for $h>0$; since $f$ is continuous, $f=0$ follows.\\
Proof of (b): One can assume $A_1=A_2=PS_+$. If $f \subset PS_+, h\in BUC, f+\gamma h=:k\in C_0$, then $M_hf=M_h(k-\gamma h) \in C_0 + C_0=C_0$ by Lemma 3.11; as in (a), $f=0$ follows. So $\gamma h \in C_0; \liminf |\gamma| >0$ and the uniform continuity of $h$ imply therefore $h\in C_0$, then $h=0$ follows as in (a) if $h\in PS_+$. --

\begin{example}
If $ \gamma $ is as in Lemma 3.11, $\gamma \in C(J,\F)$ with $\gamma (t)\ne 0$ for $t\in J, $ e.g. $g_{\omega,r}$ with $0\ne \omega \in \R, 1 < r \in \R$, and $\{0\}$ strictly $\subset A\subset PS_+(J,X)\cap BUC(J,X)$, then $\gamma A(J,X) \not\subset M(\gamma A (J,X))$, one even has $(\gamma A)\cap M(\gamma A) =\{0\}$,
though e.g. $gAP(\R,X)$ is linear, invariant and uniformly closed, so $gAP$ does also not satisfy ($\Delta$).

Further examples: $A\subset MA, $ ($\Delta$) false: Example 4.8; $A\not\subset MA, (\Delta): X+gX$ (Proposition 6.10).
\begin{proof}
If $\phi := \gamma f \in M(\gamma A)$ with $f\in A \subset BUC$, the $M_h \phi \in C_0$ by Lemma 3.11, and $M_h\phi \in \gamma A$; Lemma 3.13(b) gives $M_h \phi=0, h>0$; since $\phi$ is continuous, one gets $\phi=0$, so $f=0$.\\
``Not ($\Delta$)'' follows with Proposition 3.1, if $A$ is linear positive invariant, e.g. $A=AP$\mbox{. --} \end{proof}
\end{example}
For $\gamma, A$ as in Example 3.14 one gets then
\begin{align}
M^m(\gamma A) \cap M^n(\gamma A) = \{f \in X^J: f=0 \text{ a.e.}\}, m\ne n, m,n\in \N,\\
\text{dim} M^mgAP=\infty \text{ for all } m\in \N \nonumber.
\end{align}
Furthermore one can show
\begin{equation}
M(gAP(\R,X))=\{g'\phi + g\phi ': \phi \in AP\cap W_{\text{loc}}^{1,1}(\R,X)\} +\{f=0 \text{ a.e.} \},
\end{equation}
containing essentially only unbounded functions.\\
In the case $A=$  constant functions $X_c$ one can show
\begin{example}
For arbitrary $\gamma \in \F^J, $ the following (a), (b), (c), (d) are equivalent, with $\mathcal{N}(J,X):=\{f\in X^J: f=0 \text{ a.e.}\}$ (see also Example 4.4):\\
(a) $\gamma X \subset M(\gamma X)$ and $\gamma \notin \mathcal{N}(J,\F)$\\
(b) $M(\gamma X) \ne \mathcal{N}(J,X)$,\\
(c) $\gamma X$ is positive invariant, $\gamma \in \Lloc(J,\F), \gamma \ne \mathcal{N}(J,\F),$\\
(d) there exist $a,b \in \F$ with $a\ne 0$ and $\gamma (t)=a e^{bt}, \ t\in J.$
\end{example}
\begin{proposition}
If $\gamma \in M(o(1/|t|))\cap C(J,X)$ with $\norm{\gamma}_{S^1} <\infty$ and $\inf_J |\gamma| >0$, if further $C_0(J,X)\subset A\subset UC(J,X)$, then $\gamma A \subset M(C_0(J,X))\subset M(\gamma A).$
\begin{proof}
Since $A\subset UC$, by (3.5) to $f\in A$ and $h,\epsilon \in \R^+$ exist $n, h_j$ with \\$\norm{hM_h(\gamma f)(t) - \sum_1^n \phi_j (t)} \leq \epsilon (h+1) \norm{\gamma}_{S^1}, t\in J$, with \\$\phi_j(t)=f(t+h_j)\int_{h_{j-1}}^{h_j}\gamma(t+s)ds.$
With $C_0$ uniformly closed it is enough to show $\phi(t):=f(t+b)\int_0^{b-a} \gamma(a+t+s)ds \in C_0$ if $0<a<b$.
Now to $f\in A \subset UC$ exist $\beta $ and $ t_1>1 $ with
\begin{equation}
\norm{f(t)}\leq \beta |t| \text{ if } t_1 \leq |t|, t\in J,
\end{equation}
so $\norm{f(t+b)}\leq \beta |t+b|\leq 2\beta | t|$ if $|t|\geq t_2 =t_1+b.$
The assumptions on $\gamma$ give to $\epsilon$ and $b-a$ an $t_3\geq t_2$ with
$|\int_0^{b-a} \gamma(t+s) ds | \leq \epsilon \frac{1}{|t|}$ if $t_3 \leq |t|,t\in J$ or $|\int_0^{b-a} \gamma(a+t+s) ds| \leq \frac{\epsilon}{|t+a|} \leq \frac{\epsilon}{|t|-a}\leq \frac{2\epsilon}{|t|}$ if $|t| \geq t_3$, with $t_3 \geq 2a$.
Together one gets $\norm{\phi(t)} \leq 2\beta \  |t|  \ \frac{2\epsilon}{|t|}=4\beta \epsilon, |t| \geq t_3, t\in J$. This implies $\phi \in C_0, =\gamma \frac{1}{\gamma} C_0\subset \gamma A$ as desired.
-- \end{proof}
\end{proposition}
\begin{example}
If $A$ is as in Proposition 3.16, e.g. $A=UC(J,X)$ then with $\F=\C$
$$
g_{\omega,r} A \subset M(C_0(J,X)) \subset M(g_{\omega,r}A) \text{ if } 0 \ne \omega \in \R, 2<r\in \R
$$
\begin{proof}
with Proposition 3.16: Integration by parts yields for $h\in \R^+, |t| \geq 2h+1, t>0 \text{ resp. } t<0, $
\begin{align}
&i\omega r h (M_h(g_{\omega,r}))(t)=i \omega r \int_0^h e^{i\omega|t+s|^r}ds=\\
&\pm \left(\frac{g_{\omega,r}(t+h)}{|t+h|^{r-1}} - \frac{g_{\omega,r}(t)}{|t|^{r-1}}\right) + (r-1) \int_0^h g_{\omega, r}(t+s) |t+s|^{-r} ds,\nonumber
\end{align}
\begin{equation}
|\omega r h (M_h(g_{\omega,r}))(t)| \leq (2^r+1)\frac{1}{|t|^{r-1}}, \ |t|\geq 2h+1 \text{ . --}
\end{equation}
\end{proof}
\end{example}
If $r=2$, by (3.12) one has only $M_h(g_{\omega,2})=0(1/|t|),$ further calculations show that this cannot be improved.\\
Also, $g_{\omega ,2}UC \not\subset MC_0: M_h(g_{\omega,2}t)\notin C_0.$\\
However, one has
\begin{example}
If $0\ne \omega \in \R, 2\leq r \in \R$, and $BUC(J,X) \subset A\subset UC(J,X),$ then $g_{\omega,r} A \subset M(g_{\omega,r} BUC(J,X))\subset M(g_{\omega,r}A)$.
\end{example}
\textit{Proof,} $r>2$: Example 3.17.\\
Proof, $r=2$, due to Bolis Basit, May 2013, unpublished:
For $h\in \R^+, t\in J, f\in UC(J,X)$ one has
\begin{align*} hM_h(g_{\omega,2}f)(t)=g_{\omega,2}(t)\int_0^h g_{\omega,2} (s) \gamma_{2\omega s} (t) (f(t+s)-f(t)) ds +\\
g_{\omega,2} (t) f(t) \int_0^h g_{\omega,2}(s) \gamma_{2\omega t}(s) ds =:I_1 + I_2.
\end{align*}
To $f\in UC$ exists $n_1 \in \N $ with $|\Delta_u f| \leq 1 $ if $0\leq u \leq 1/n_1$,  implying $|\Delta_s f | \leq $\\ $(h+1)n_1$ on $J$ if $0\leq s\leq h$.
If $|\Delta_u f| \leq \epsilon $ on $J$ if $0\leq u\leq \delta$, then $\norm{(\Delta_sf)(t+u) - (\Delta_s f)(t)}$
$\leq \norm{f(s+t+u)-f(s+t)}+\norm{f(t+u)-f(t)}\leq 2\epsilon$, $t\in J, 0\leq s \leq h$.
Similarly, $\gamma_{2\omega s}(t)$ is bounded and uniformly continuous in $t$ on $J$, uniformly in $0\leq s\leq h$. The above implies the same for $\gamma_{2\omega s}(t) \Delta_s f(t), $ so $I_1 \in BUC(J,X)$.

\begin{align}
&2i \omega I_2=f(t)g_{\omega,2} (t) \frac{1}{t}\int_0^h g_{\omega,2}(s) 2 i \omega t e^{i \omega t s } ds=\\
&=g_{\omega,2} (t) \frac{f(t)}{t} [g_{\omega,2}(h) \gamma_{2\omega h }(t)-1] - 2i \omega \int_0^h sg_{\omega,2}(s) \ \gamma_{2 \omega s}(t)ds.\nonumber
\end{align}
Now (3.10) and some simple calculations give $f(t)/t$, $[\dots]$ and $\int_0^h sg_{\omega,2}(s) \gamma_{2\omega s}(t) ds$ bounded and uniformly continuous on $|t| >1,$ so $I_2\in BUC(J,X)$, also if $J=(\alpha,\infty).$
--

\begin{example}
If $1<r\in \R, \omega \in \R,$ then $ g_{\omega,r}UC(J,X)\subset M(g_{\omega,r}UC(J,X)),$ but not $\subset M(g_{\omega,r}BUC(J,X))$ if $r<2$.\end{example}
\noindent This follows with a refinement of Bolis Basit's proof for $r=2$ (Example 3.18), we omit the somewhat lengthy details.
\begin{example}
If $\gamma=g_{\omega,1}, \ \omega \in \R,$ or $\gamma=\gamma_{\omega}, \ \omega \in \C$, also $\gamma UC(J,X)\subset M(\gamma UC(J,X)),$ but not $\subset M(\gamma BUC(J,X)) $ if $\omega\ne 0$.
\end{example}
\abstand
\section{$\Delta$-classes and the ($\Delta$)-condition}
The ($\Delta$)-condition (2.13) has been introduced in \cite[Def. 1.4 p. 7]{lit5}, as indicated in the introduction it has been useful in several instances.\\
By Proposition 3.1, for linear positive invariant $A,A\subset MA$ is necessary for ($\Delta$). Sufficient conditions are given by
\begin{proposition}
If $A$ is linear, positive invariant, uniformly closed and $A\subset UC(J,X)$, then $A$ satisfies ($\Delta$).
\begin{proof}
\cite[Theorem 2.4 p. 428]{lit9}.
-- \end{proof}
\end{proposition}
Examples are $C_0, AP, $ uniformly continuous almost automorphic or Levitan almost periodic functions, Eberlein weakly almost periodic functions $EAP,BUC,UC$ \cite[p. 8 (1.2a)-(1.2f)]{lit6}, \cite[Ex.s 3.2, 3.3]{lit8}.\\
Direct proofs of ($\Delta$) have been given for $BC(J,X)$ \cite[Ex.3.5 p. 1013]{lit8}, weighted $L^p$ spaces \cite[Prop. 3.4 p. 1012/1013]{lit8}, various spaces of ergodic functions \cite[Prop.s 3.8 and 3.10 p. 1014/1015]{lit8}, classes of functions satisfying a Lipschitz condition or a growth condition $O(w)$ with suitable $w$ \cite[Ex.s 3.7, 3.13 p. 1013/1015]{lit8}.\\
The $g_{\omega,r}AP$ or more generally $\gamma A$ of Example 3.14 do not satisfy ($\Delta$) if they are positive invariant, since they do not satisfy $A\subset MA$ (Proposition 3.1).\\

To get ($\Delta$) for $\gamma BUC$ or $\gamma UC$, we need a category result:
\begin{lemma} Assume $\beta \in \R^+$ and $A, A_n \subset \Lloc (J,X), n\in \N, $ with $A=\cup_1^{\infty} A_n, \ f\in \Lloc (J,X),$ with:
\begin{align}
&\Delta_sf \in A \text{ for } 0 \leq s\leq \beta,\\
&\text{if } s_m \in [0,\beta], m\in\N, \text{ with } s_m \to r \text{ and } \Delta_{s_m} f \in A_n \text{ with fixed } n,\\
&m\in \N, \text{ then } \Delta_rf \in \text{ this } A_n \text{ too},\nonumber \\
&\text{if with fixed } n\in\N \text{ and } u, \rho \text{ with  } 0\leq u< u+\rho \leq \beta \text{ one has } \\
&\Delta_s f \in A_n \text{ for } u\leq s\leq u + \rho, \text{ then } \frac{1}{\rho}\int_u^{u+\rho} \Delta_s f ds \in A_n \  [\text{ resp. } A]\nonumber
\end{align}
Then there exist $v, \delta$ with $0 \leq v < v+ \delta \leq \beta$ and $m \in \N$ such that $(M_h f)_v - f \in A_m$ [resp. $A$] for $0<h \leq \delta$.
\begin{proof}
\cite[Lemma 3.6, p. 1013]{lit8}, \cite[Lemma 4.12, p. 33/34]{lit6}.
-- \end{proof}
\end{lemma}
Here in (4.3) the integral is meant as the Bochner integral of the $X$-valued function $s\to f(t+s)-f(t)$, with fixed $t$.\\

Let us assume now $A=\gamma UC(J,X)$ with
\begin{align}
&\gamma \in C(J,\F), \gamma(t)\ne 0 \text{ for } t\in J,\\
&f \in \Lloc (J,X), \text{ with } \Delta_s f\in A \text{ for } 0\leq s\leq \beta, \beta \in \R^+.
\end{align}
By Lemma 2.1, (4.5) implies
\begin{equation}
f\in C(J,X)
\end{equation}
With (4.4), (4.5), (4.6) one gets
\begin{align}
&(\Delta_s f)(t)=\gamma(t) \Phi (s,t) \text{ with unique  } \Phi, \Phi \in C([0,\beta] \times J, X),\\
& \Phi(s, \cdot)\in UC(J,X), \  0\leq s\leq \beta. \nonumber
\end{align}
With $\epsilon \in \R^+, $ fixed in the following, and $n\in \N$, define
\begin{align}
&A_n^{\epsilon}:=\{\gamma \phi: \phi \in C(J,X) \text{ with } \norm{\phi(t+h)-\phi(t)}\leq \epsilon \text{ if } 0\leq h\leq \frac{1}{n}, t\in J\},\\
&A^{\epsilon}:=\cup_{n=1}^{\infty} A_n^{\epsilon}
\end{align}
Lemma 4.2 can be applied with $A=A^{\epsilon}, A_n=A_n^{\epsilon}:$\\
(4.1) follows with (4.5)-(4.8).\\
(4.2) Continuity (4.6) of $f$ implies $s$-continuity of $\Phi$ of (4.7), so $\Phi(s_m, \cdot) \in A_n^{\epsilon}$ for $m\in\N $ implies $\Phi(r,\cdot) \in \text{ same } A_n^{\epsilon}$.\\
(4.3), case $A_n$: By (4.7) the integral  $I$ in (4.3) exists, $I=\gamma(t) \frac{1}{\rho} \int_u^{u+\rho}\Phi(s,t)ds$, with $\norm{\Phi(s,t+h)-\Phi(s,t)}\leq \frac{1}{n} $ if $0\leq h\leq \frac{1}{n}$ and $t \in J$, any $s\in [u,u+\rho]$, this gives $I\in A_n^{\epsilon}$.\\
So there exist $m\in\N$ and $v, \delta $ with $0\leq v<v+\delta\leq \beta$ with
\begin{equation}
(M_hf)_v-f\in A_m^{\epsilon},\  0 <h\leq \delta.
\end{equation}
Next we assume, with suitable $\gamma$, that
\begin{equation}
\gamma UC(J,X) \subset M(\gamma UC(J,X)).
\end{equation}
Then $-M_h(\Delta_v f) \subset \gamma UC(J,X) \subset A^{\epsilon}$ of (4.9), so with (4.10)
\begin{align}
&M_h(f)-f=M_h(f_v)-f-M_h(\Delta_v f)=(M_h(f))_v-f-M_h(\Delta_v f)\in \\
&A_m^{\epsilon} + A^{\epsilon} \subset A^{\epsilon}+A^{\epsilon}\subset A^{2\epsilon}, \  0< h \leq \delta. \nonumber
\end{align}
Now for any $f\in \Lloc(J,X)$ and $h\in \R^+$ one has
\begin{equation}
-f+M_{2h}f = -f+M_h f - \frac{1}{2}M_h(\Delta_hf).
\end{equation}
Since $\Delta_s f \in A $ by (4.5), with (4.11) and (4.9) one has $M_h(\Delta_sf) \in A \subset A^{\eta}, h,s,\eta \in \R^+$; choosing $\epsilon_j \in \R^+ $ with  $\sum_1^{\infty} \epsilon_j <\epsilon$, induction yields with (4.12)
\begin{align}
& M_{2^n h}f-f\in A^{2\epsilon} + \sum_1^{n} A^{\epsilon_j} \subset A^{2\epsilon} + A^{\sum_1^n \epsilon_j} \subset A^{2\epsilon} + A^{\epsilon} \subset A^{3\epsilon},\\
& 0 < h\leq \delta, n\in \N . \nonumber
\end{align}
So $ f-M_hf \in -A^{3\epsilon} =A^{3\epsilon}$ for any $h\in \R^+$. \\
$\epsilon$ being arbitrary, one gets $f-M_hf \in \cap_1^{\infty} A^{1/m}$; with (4.8), (4.7) this implies $f-M_hf\in \gamma UC =A, h\in \R^+, $ i.e. ($\Delta$) holds for $\gamma UC(J,X).$\\
We have now shown (a) of the following
\begin{theorem}
Assume $ \gamma \in C(J,\F)$ with $\gamma(t)\ne 0$ for $t\in J$. Then\\
(a) $\gamma UC(J,X) \subset M(\gamma UC(J,X)) $ implies ($\Delta$) for $\gamma UC(J,X)$.\\
(b) $\gamma BUC(J,X) \subset M(\gamma BUC(J,X))$ implies ($\Delta$) for $\gamma BUC(J,X)$.
\begin{proof}of (b): (Since we do not assume (4.11), (a) cannot be used): If in the above proof of (a) in (4.5) the $A$ is $\gamma BUC$ and (4.11) is replaced by $\gamma BUC(J,X)\subset M(\gamma BUC(J,X))$, if further in the definition (4.8) of $A_n^{\epsilon}$ one adds ``$\norm{\phi(t)}\leq n$ for $t\in J$'', with corresponding $A^{\epsilon}$, one gets again $f-M_h f\in \cap_1^{\infty} A^{1/m}$, which now means $f-M_h f\in \gamma BUC(J,X)$.
-- \end{proof}
\end{theorem}
Examples for $\gamma$  with $\gamma BUC \subset M(\gamma BUC)$ are given in Proposition 3.12 and Lemma 3.10, for $\gamma UC \subset M(\gamma BUC)$  in Proposition 3.16/Examples 3.17/3.18/3.19/ 3.20, especially
\begin{align}
& g_{\omega,r} BUC \subset M(g_{\omega, r} BUC) \text{ and } g_{\omega, r} UC \subset M(g_{\omega, r} UC)  \text{  if } \omega \in \R, 1\leq r\in \R, \\
& \text{ so these }  g_{\omega, r} BUC  \text{ and    } g_{\omega,r} UC  \text{ satisfy }  (\Delta),  \text{ all }  J  \text{ and } X. \nonumber
\end{align}
\begin{example}
$A=\gamma X$ with $\gamma \in C^1(J,X), \gamma \notin \{\rho e^{\beta t} : \beta, \rho \in \F\}$ is linear, uniformly closed and with $\Delta A = X_c, $ so $A$ satisfies ($\Delta$), but $MA=\{f=0  \text{ a.e.} \}$, so $ A\not\subset MA$. This also shows that the assumption ``positive invariant'' is essential in Proposition 3.1. See also Example 3.15.
\begin{proof}
$f \in \Delta A$ of (2.12) means $ f(t+h)-f(t)=\gamma (t) a(h), t\in J, h\in \R^+, a(h)\in X$; then $f\in C^1(J,X)$ by Lemma 2.1, with $\gamma \not\equiv 0$ also $a\in C^1(\R^+, X), $ so $f'(t+h)=\gamma (t) a'(h). $ With $f'\in C^1$ and $\gamma \not\equiv 0, b:=a'(0+)$ exists $\in X$,
$f'(t)=\gamma (t) b, t\in J$. With this one gets
$\gamma  (t) a(h)=f(t+h)-f(t)=b\int_t^{t+h} \gamma (s) ds$, then $\gamma(t) a'(h)=b\gamma (t+h)$; as above, this implies $a'\in C^1 $ and $\gamma (t) a''(h)=b\gamma'(t+h), c:=a''(0+)$ exists $\in X$ and $\gamma '(t)b=\gamma(t) c, t\in J$. If $b \ne 0$, Hahn-Banach gives $\beta \in\F$ with $\gamma'(t)=\beta \gamma(t)$ or $\gamma (t)=\rho e^{\beta t}$ against the assumption, so $b=0$. Then $f'(t)=\gamma(t) b = 0, f\in X_c,$ indeed $\Delta(\gamma X)=X_c$. Then by definition (2.13), the $X$ satisfies ($\Delta$). $MA=\{f=0 \text{ a.e.}\}$ follows with Example 3.15.
-- \end{proof}
\end{example}
Concerning general $\Delta$-spaces $A$, the following is obvious
\begin{equation}
\text{ If } A\subset X^J \text{ is linear and positive invariant, then } A\subset \Delta A \subset \Delta ^2 A \subset \dots .
\end{equation}
As in Proposition 3.3 for $MA$, one has (with $\Delta_h \Delta_k f = \Delta_k \Delta_h f, \Delta_h M_kf =M_k(\Delta_hf))$
\begin{proposition}
If $A \subset \Lloc (J,X)$ is linear resp. positive invariant resp. invariant resp. uniformly closed resp. satisfies ($\Delta_1$) resp. satisfies ($\Delta$) resp. satisfies $A\subset MA$, the same holds for $\Delta A$.
\end{proposition}
\begin{proposition}
If $A$ satisfies ($\Delta_1$) and $\R^+ A \subset A$, then
\begin{equation}
\Delta A \subset A+X_c + PA,
\end{equation}
with ($t_0 $ fixed $\in J$)
\begin{equation}
PA:=\{Pf: f\in A\cap \Lloc (J,X)\}, (Pf)(t):=\int_{t_0}^t f(s) ds, \  t\in J.
\end{equation}
\begin{proof}(see [7, Prop. 3.3])
If $f\in \Delta A$, the definitions give $f\in \Lloc$ and $f-M_hf=:u\in A$ for some $h\in \R^+$, so $hf(t)=hu+\int_t^{t+h} f ds = hu+\int_{t_0 + h }^{t+h} f ds - \int_{t_0+h}^t f ds= hu+ \int_{t_0}^t \Delta_h f ds + \int_{t_0}^{t_0+h} fds.$
-- \end{proof}
\end{proposition}
Proposition 4.6 implies (for condition ($\Delta P$) see §8)
\begin{align}
&\text{If } A \text{ is linear positive invariant with } (\Delta_1) \text{ and } A\subset MA, \text{ then } \\
& \Delta A= A+X_c + PA. \nonumber \\
&\text{If } A \text{ is as in (4.19) and satisfies } (\Delta P), \text{ then } A\text{ satisfies } (\Delta).
\end{align}
\begin{example}
$\Delta(BUC(J,X))=UC(J,X)=BUC(J,X)+P(BUC(J,X)).$
\begin{proof}
\cite[Prop. 4.1 p. 1010]{lit8} for first ``$\subset$''; ``$\supset$'': $|\Delta_s f|\leq (1+h)n_1$ in Proof of Ex. 3.18; (4.19).
-- \end{proof}
\end{example}

\begin{example}
There exist $A \subset BC(\R,X)$ which are linear, invariant, uniformly closed with $A\subset MA$, but they do not satisfy ($\Delta$):\\
$BEM(J,X), BEM_0(J,X), CEM(J,X),$ $CEM_0(J,X)$ with $J=\R_+$ or $\R$ and any Banachspace $X$ are all
\begin{align} &\text{linear, positive invariant (invariant if }J=\R \text{), uniformly closed,}\\
 &\text{with } A\subset MA,\nonumber
\end{align}
but \emph{none} of them satisfies ($\Delta$).
\end{example}
Here $f\in EM(J,X)$ or $f$ is \textit{Maak-ergodic} means to $f\in X^J$ exists $a\in X$ so that to any $\epsilon \in \R^+$ exist $n\in \N$ and $s_1, \dots, s_n \in J$ so that
\begin{equation}
\norm{\frac{1}{n} \sum_1^n f(s_j + t) - a} \leq \epsilon \text{ for all  } t\in J.
\end{equation}
The $a$ is then unique, $m(f):=a$  (\cite[p. 373, Mittelwertsatz]{lit30});\\
$CEM(J,X):=EM(J,X)\cap C(J,X), BEM:=EM\cap BC, EM_0:=\{f\in EM: m(f)=0\}$. $E(J,X):=\{\text{uniformly ergodic } f: J\to X\}$ (see e.g. \cite[p.117/118]{lit4}), $BE:=E\cap BC.$\\
$BEM$ has not ($\Delta$): Assume $f\in BE \subset \Lloc, $ then $\Delta_h f \in BEM$ for $h \in \R^+$ by Remark 4.10, so assuming ($\Delta$) for BEM one has $f-M_hf \in BEM$; now $BE\subset MBE$ (Fubini), so $M_h f \in BE, f$ bounded implies $M_h f \in BUC; $ so $ M_h f \in BE\cap BUC, = EM \cap BUC$ by Remark 4.9, $\subset BEM;$ this implies $f\in BEM, $ and so $BE\subset BEM$. But one can construct $f\in BE_0$ with $f\notin EM$, a contradiction. So $BEM$ has not ($\Delta$).\\
If $f\in BE_0,$ the same proof shows that $BEM_0$ has not ($\Delta$). ($\Delta$) for $CEM$ would with ($\Delta$) for $BC$ (Lemma 2.1 and \cite[Prop.1.1]{lit7}) imply ($\Delta$) for $BEM$, similarly for $CEM_0$. --\\

Construction of an $f_0\in BE_0, f_0 \notin EM:$ With $\phi(t):=e^{-t^2}$ and $v_1:=2$ define recursively $v_{n+1}:=v_n + \phi(v_n), n\in \N;$ by a contradiction argument, $v_n \to \infty;$ with $w_n:=v_n + (1/v_n)\phi(v_n), a_n:=\frac{1}{8v_n} \phi(v_n)$ define $f_n=1$ on $[v_n, w_n]$, $=0 $ outside $(v_n-a_n, w_n + a_n)$, linear else.
Then $\text{supp}f_m \cap \text{supp}f_n =\emptyset$ if $m\ne n$, so $f_0:=\sum_1^{\infty} f_n$ is well defined,
$\in C(\R, [0,1]).$ With $0<\int_{v_n}^{v_n + k} f_0 ds < \frac{5}{4v_n}(v_{n+k}-v_n) $ one gets $f_0 \in BE_0$.
Assuming $f_0\in EM, $ one has $m(f)=0$ since on $E\cap EM$ the mean on $E$ coincides with the $m$ for $EM$; so to $\epsilon=1/2$ there exist $n$ and $s_1\leq s_2 \leq \dots \leq s_n$ (multiplicities can appear) for which (4.22) holds with $a=0$.
If $2\delta := \min \{s_{j+1}-s_j: 1\leq j\leq n-1 \text{ and } s_{j+1}-s_j >0\}>0$ and if there is $n_m \in \N$ with $f(s_j + t) =1 $ for all $t\in U_m:=[v_{n_m}, w_{n_m}]-s_m$ and $1\leq j\leq m <n$ and additionally $\phi(v_{n_m})\leq \delta$ and $v_{n_m} \geq \delta^{-2}, $ then one can find $n_{m+1}>n_m$ so that $U_{m+1}\subset U_m $ and $U_{m+1}$ has the same properties as $U_m$, but for $s_j$ with $1\leq j\leq m+1$. With induction one gets for $m=n$ a contradiction with the assumed (4.22), $a=0, \epsilon=1/2$.\\
For $n_{m+1}$, assuming $s_{m+1}>s_m$, choose $p$ maximal with $v_p\leq v_{n_m} + s_{m+1}-s_m, $ then $n_{m+1}:=p+1$
does it, using $e^{-(t+\delta)^2} \leq \frac{1}{2t}e^{-t^2}4$ if $t\geq \delta^{-2}$. --
\begin{remark}
With the notation after Example 4.8 one has $BEM\subset EM\cap ML^{\infty} \subset E$ (integrate (4.22) from $s$ to $s+T$), with the $f_0$ of Example 4.8 the second ``$\subset$'' is strict. One has however
$EM\cap UC=EM\cap BUC =E\cap BUC=E \cap UC$, these four spaces satisfy ($\Delta$) (Proposition 4.1); the same holds if $EM$ resp. $E$ are replaced by $EM_0$ resp. $E_0$.
\end{remark}
\begin{remark}
If $f\in X^J$ with $\norm{f}_{\infty}< \infty, h\in \R^+$, then $ \Delta_h f \in EM_0 (J,X), J\subset \R^+$ or $= \R:$
$\frac{1}{n} \sum_1^n(\Delta_h f)(jh+t)=\frac{1}{n}(f((n+1)h+t)-f(t))\to 0$ uniformly in $ t\in J$ (B. Basit and A. Pryde, Analysis Paper 101 Feb. 1996, Monash Univ.). \\

More examples of $A$ with ($\Delta$) can be found in \cite[Ex./Prop. 4.7-4.11,4.14-4.17]{lit6}, \cite[Ex. 3.14-3.19]{lit7}, \cite[Ex. 3.2-3.5, 3.8, 3.10, 3.13]{lit8}, and in §§ 6 and 8 below.
\end{remark}
\abstand
\section{Vector sums}
In the following we discuss $A=U+\gamma V:=\{u+\gamma v: u\in U, v\in V\}$, mostly with (36 pairs $(U,V)$, general $J, X$)
\begin{align}
&U,V \in \mathcal{U}=\{ X_c, C_0(J,X), X_c + C_0 (J,X), AP(J, X), BUC(J,X), UC(J,X)\},\\ &\gamma \in \F^J.\nonumber
\end{align}
\textbf{1.} For $U,V$ of (5.1), the $U+\gamma V$ is always linear, it is $\subset \Lloc (J,X)$ resp. $C(J,X)$ if $\gamma $ is.\\

$U+\gamma V$ is a direct sum, $U\cap \gamma V=\{0\}$, in the following cases:
\begin{enumerate}
\item[(a)] $X+\gamma V, AP+\gamma V, $ any $V\in \mathcal{U}, \gamma \in MC_0(J,X), \norm{\gamma}_{S^1}<\infty$, and additionally
$\liminf_{|t|\to\infty} |\gamma| >0 $ if $V=UC$ (Lemma 3.13(a)).
\item[(b)] $C_0 + \gamma X, X_c + C_0 + \gamma X,  C_0 + \gamma AP, $ if $\liminf_{|t|\to\infty}|\gamma| >0$.
\item[(c)] If $U=BUC$ resp. $UC, \gamma \not\equiv 0, $ then $U+\gamma X$ is direct if and only if $\gamma \notin BUC(J,X)$ resp. $\notin UC(J,X)$.
\item[(d)] $U+ \gamma AP, U=BUC $ or $UC$, $\liminf_{|t|\to\infty}|\gamma|>0, \gamma $ as in Lemma 5.1. \\
(Lemma 5.1 below gives $\gamma f\in C_0, \liminf |\gamma| >0$ then $f\in C_0$, so $f=0$ if $f\in AP$.)
\end{enumerate}
$\gamma=g_{\omega, r}$ with $r>1, 0 \ne \omega \in \R$ fulfill all the assumptions in (a)-(d) (Lemma 3.10, after Lemma 5.1).\\
The ``remaining'' 16 cases all contain $C_0 + \gamma C_0$, so for these $U+\gamma V$ is not direct if $0\not\equiv \gamma \in C(J,X)$.
\begin{lemma}
If $\gamma \in \F^J$ satisfys $O_1$ of (5.3) and $\norm{\gamma}_{\infty} <\infty$, then
\begin{equation}
UC(J,X)\cap (\gamma UC(J,X)) \subset C_0(J,X).
\end{equation}\end{lemma}
Here $\gamma$ \textit{satisfies the (oscillatory)  condition } $O_1$ means
\begin{align}
&\gamma \in \F^J, \text{ to each sequence } (t_n)_{n\in \N} \text{ from } J \text{ with } |t_n| \to \infty \text{ there} \\
&\text{exists } \rho_0 \in \R^+ \text{ so that to each } \delta \in \R^+ \text{ there exists} \nonumber \\
& m_{\delta} \in \N \text{ and } \nonumber
s_{+}, s_{-} \in [t_{m_{\delta}}, t_{m_{\delta}} + \delta] \text{with } |\gamma (s_+) - \gamma (s_-) | \geq \rho_0.
\nonumber
\end{align}
Lemma 5.1 is applicable with $\gamma = \phi g_{\omega, r}$ of (2.3), $0\ne \omega \in \R, 1<r\in \R$, $\phi \in BU(J,X)$ with
$\inf_J |\phi| >0$ (Lemma 3.10, (5.4)); $\gamma=\sin (\omega |t|^r)$ satisfies $O_1$, same $\omega, r$. See also Proposition 5.5.
\begin{equation}
\text{If } \gamma \in MC_0(J,\F) \text{ and } \liminf_{|t|\to \infty}|\gamma | >0 \text{ then } \gamma \text{ satisfies } O_1.
\end{equation}
So $\gamma \in MC_0$ is stronger than $O_1$ (we omit the proof of (5.4));\\
$\gamma = 2 + g_{1,2}$ satisfies $O_1$, but $\gamma \notin MC_0$.\\
If in Lemma 5.1 additionally $\gamma \in C(J,\F)$ with $\inf_J |\gamma| >0,$ then one has equality in (5.2).\\
Without $O_1$ resp. $\norm{\gamma}_{\infty}<\infty$ Lemma 5.1 becomes false: $\gamma =\frac{1}{1+|t|}$ with
$u=1$ (or $\gamma\equiv 1$) resp. $\gamma =\sqrt{1 + |t|} g_{1,2} $ ($\in MC_0$ by Lemma 3.10, (5.4), $1 \in UC\cap \gamma UC$).\\

\textit{Proof of Lemma 5.1 by contradiction:} If $u,\gamma v \in UC$ with $u=\gamma v$, but $v\notin C_0,$ there exist $t_m\in J$ with
$|t_m|\to\infty$ and $\epsilon_0 \in \R^+$ with $|v(t_m)|\geq 4 \epsilon_0, m \in \N$. By $O_1$, there exists $\rho_0 \in \R^+$
so that to any $\delta\in\R^+$ there exist
\begin{equation}
m_{\delta} \in \N, s_+, s_- \in [t_{m_{\delta}}, t_{m_{\delta}} + \delta ] \text{ with } |\gamma(s_+) - \gamma (s_-)| \geq \rho_0;
\end{equation}
furthermore, to $\epsilon:=\min \{\epsilon_0, \rho_0 \epsilon_0\}$ there is $\delta_{\epsilon} \in \R^+$ with
$\norm{u(t)-u(s)} \leq \epsilon$ and $\norm{\gamma}_{\infty} \norm{v(t)-v(s)}\leq \epsilon$ if $s,t \in J$ with $|s-t|\leq \delta_{\epsilon}$;
with the above $\delta=$ this $\delta_{\epsilon}$, one has (5.5). Then
\begin{align*}
\rho_0 \epsilon_0 \geq &\epsilon \geq \norm{u(s_+)-u(s_-)}  =\norm{\gamma (s_+ ) v(s_+)- \gamma (s_-) v(s_-)} \geq \\
& \norm{(\gamma (s_+)-\gamma (s_-)v(s_+))} - \norm{\gamma (s_-)(v(s_+)-v(s_-))} \geq \rho_0 \norm{v(s_+)} - \epsilon \geq \\
&\rho_0(\norm{t_{m_{\delta}}} - \norm{v(s_+)-v(t_{m_{\delta}})}) - \rho_0 \epsilon_0 \geq \rho_0 (4 \epsilon_0 - \epsilon_0 -\epsilon_0)=2\rho_0 \epsilon_0 >\rho_0\epsilon_0 \text{. --}
\end{align*}
\begin{corollary}
If $\gamma $ is as in Lemma 5.1 and $\liminf_{t \to \infty}>0$, then (see (3.6))
\begin{equation}
(\gamma(PS_+(J,X)\cap UC(J,X))) \cap UC(J,X)=\{0\}
\end{equation}
\normalfont So with (5.4), $\gamma \in MC_0$, $\norm{\gamma}_{\infty}<\infty$ and $\liminf_{|t|\to\infty} |\gamma | >0$ give (5.6);\\
without $\liminf |\gamma| >0$ this is false. (See also (3.7)).

\textit{Proof of Cor. 5.2:} If $f\in (\gamma ( PS_+\cap UC))\cap UC, f\in C_0$ by Lemma 5.1. With $f=\gamma \phi, \phi \in PS_+, $ and $|\gamma | \geq \delta_0>0$ on some $[n,\infty)$
with $\liminf |\gamma| >0$, one gets $\phi(t)\to 0$ as $t\to \infty$: the definition (3.6) of $PS_+$ gives $\phi=0$.--
\end{corollary}

\textbf{2.} $U+\gamma V$ is positive invariant if $\gamma V$ is, any $U\in \mathcal{U}$.
\begin{lemma}
If $\gamma \in \F^J$ with $\gamma(t)\ne 0$ for $t \in J, A$ is [positive] invariant $\subset X^J$ with $(y \circ u)v\in A$ if $u,v \in A, y \in \textnormal{ dual } X'$, and there exists $x_0 \in A, 0\ne x_0\in X_c$, then $\gamma A$ is [positive] invaviant if and only if
\begin{equation}
\text{there is } \epsilon_0\in \R^+ \text{ with } \gamma_a x_0 /\gamma \in A \text{ for all } a\in [-\epsilon_0, \epsilon_0] \ [\text{resp. } a\in (0, \epsilon_0)].
\end{equation}
\begin{proof}
The necessity of (5.7) is obvious. Conversely, (5.7) implies $\gamma_{2a}x_0/\gamma_a \in A, |a|\leq \epsilon_0;$ with $y \in X'$ with $y(x_0)=1$ one gets
$$ A \ni y(\gamma_{2a} x_0/\gamma_a)\gamma_a x_0/\gamma  = \gamma_{2a}x_0/\gamma, $$
then $\gamma _{2^na} x_0 /\gamma \in A$.
-- \end{proof}
\end{lemma}
Examples of $A$ with $(y\circ u) v \in A$: $AP(J,X), BUC(J,X)$, here $\gamma_a x_0/\gamma \in A$ can be replaced by $\gamma_a/\gamma \in AP(J,\F)$ resp. $BUC(J, \F)$, similarly for $A=X_c $ or $A=X_c + C_0 (J,X)$.\\

$g_{\omega, r} V$ is positive invariant in the following cases ($\omega \ne 0$):
\begin{enumerate}
\item[(a)] $V=X$: $r=1$ and $J\subset \R_+$; if $J\not\subset \R^+$ or $0<r\ne 1$, $g_{\omega, r} X$ is not positive invariant
$((t+h)^r-t^r=r(t+\delta h)^{r-1}, 0<\delta<1).$
\item[(b)] $\gamma (t)=e^{wt}, $ any $w\in \F, $ any $V\in \mathcal{U}$, any $J$.
\item[(c)] $V=C_0(J,X): \gamma =g_{\omega,r}, \omega \in \R, 0\leq r\in \R, $ any $J;$ for no $r<0, \omega \ne 0.$
\item[(d)] $V=AP(J,X): r=1, J\subset \R^+; r=2, $ any $J$. If $r\in (0,1) \cup (1,2) \cup (2,\infty), \omega \ne 0, g_{\omega,r}AP$ is not positive invariant.
\item[(e)] $V=BUC(J,X): 0\leq r\leq 2 \ (g_{\omega,r}(t+a)=g_{\omega,r}(t)e^{i\omega ( |t+a|^r-|t|^r)}, |t+a|^r- |t|^r$ has bounded derivative); $g_{\omega,r} BUC$ is not positive invariant if $r>2$.
\item[(f)]$V=UC(J,X): 0\leq r<1; g_{\omega,r}UC$ is not positive invariant if $r\geq 1$. \\
\end{enumerate}

\textbf{3.} $U+\gamma V$ uniformly closed: Here one cannot argue as in 1. or 2. . For $\gamma \equiv 1, U=AP$ and $V=\{$Eberlein weakly almost periodic nulfunctions$\}$, completeness of $U+V$ has been shown e.g. by Porada \cite{lit27} with what we call a Porada inequality, namely
\begin{equation}
	\norm{u}_{\infty} \leq \norm{u+v}_{\infty}, \  u \in AP(\R,\C), v\in EAP_0 (\R, \C)
\end{equation}
(see also \cite[p. 427 after (1.7), and Prop. 1.2]{lit9}, \cite[p. 51 Prop. 7.13]{lit6}).\\
For general $U+\gamma V$ (5.8) no longer holds, e.g. if $U\cap \gamma V \ne \{0\}$ (see after (5.4)).\\
For $U,V$ and $\gamma $ we are interested in however a weakened form of (5.8) still is true:
\begin{proposition}[Porada inequality] If $\gamma \in \F^J$ satisfies $\norm{\gamma}_{\infty}<\infty$
and $O_2 $ of (5.10), then
\begin{equation}
\limsup_{|t|\to\infty} |u| \leq \norm{u+\gamma v}_{\infty}, \  u, v \in UC(J,X)
\end{equation}
If only $\norm{\gamma}_{S^1} <\infty$ and $\gamma \in MC_0(J,\F)$ (instead of $O_2$), then (5.9) holds for $u\in UC$ and $v\in BUC$.
\end{proposition}
Here $\gamma$ \textit{satisfies the} (oscillatory) \textit{condition} $O_2$ means
\begin{align}
&\gamma \in \F^J, \text{ to each } (t_m)_{m\in\N} \text{ from } J \text{ with } |t_m| \to \infty \text{ there exists } \\
&\rho_0 \in \R^+ \text{ so that to each } n\in \N \text{ and } \delta \in \R^+ \text{ there exist } m_{\delta} \text{ with } \nonumber \\
&n\leq m_{\delta} \in \N \text{ and } s_+, s_- \in [t_{m_{\delta}}, t_{m_{\delta}} + \delta ] \text{ with } \nonumber \\ &|\gamma (s_+)| \geq \rho_0 \text{ and }  |\gamma (s_+) + \gamma (s_-) |\leq \delta.\nonumber
\end{align}
\begin{examples} $\gamma$ with $\norm{\gamma}_{\infty}<\infty $ and $O_2$ are again $\gamma=\phi g_{\omega, r}, 0\ne \omega \in \R, 1 <r\in \R, \phi \in BUC(J,\C)$ with $\inf_J |\phi|>0,  \F=\C; \gamma = \sin ( \omega |t|^r)$ satisfies $O_2$, same $\omega,r$.
\end{examples}
$O_2$ implies $O_1$ of (5.3): If $|\gamma (s_+)|\geq \rho_0$ and $|\gamma(s_+)+\gamma(s_-)|\leq \delta\leq \rho_0$, then $|\gamma(s_+)-\gamma(s_-)| = |2\gamma(s_+)-(\gamma(s_+)+\gamma(s_-))|\geq 2\rho_0-\delta\geq \rho_0.$
\begin{remark}
If $\gamma$ satisfies $O_2, \gamma (t) \ne 0 $ for $t\in J$ and $\norm{\gamma}_{\infty}<\infty$, then $1/\gamma$ also satisfies $O_2$.
\end{remark}
Proposition 5.4 becomes false without $\norm{\gamma}_{\infty}<\infty (\gamma=\sqrt{1+|t|}g_{1,2}$ as after (5.4), Proposition 5.8), or with $O_2$ replaced by $O_1$ ($\gamma=2+g_{1,2}, u=2, v=-1$), or with $\rho_0=0$ in $O_2 (\gamma = 1/(1+|t|), u=1,v=-(1+|t|))$.
So $O_1$ does not imply $O_2$.\\
Also, contrary to (5.4), even with $\gamma \in BC(J,\C)$ with $\inf_J |\gamma| >0,$ one can construct examples showing that
\begin{equation}
\text{neither } \gamma \in MC_0 \text{ implies } O_2 \text{ for } \gamma \text{ nor } O_2 \text{ for } \gamma \text{ implies } \gamma \in MC_0.
\end{equation}
If in $O_2$ one can get even $\gamma(s_+)+\gamma(s_-)=0$, then one can omit the $\rho_0$ and $n$ (we omit the proof):\\
\textbf{Proposition 5.4$_0$}:
If $\gamma $ satisfies $O_0$ and $\norm{\gamma}_{\infty}<\infty$, (5.9) holds.\\

$O_0$: To each $(t_m)\subset J$ with $|t_m|\to\infty$ and $\delta \in \R^+$ exist $m_{\delta}$ and $s_+,s_-\in$ \\
		${} \quad \quad \quad [t_{m_{\delta}},t_{m_{\delta}} + \delta]$ with $\gamma(s_+)+\gamma(s_-)=0.$\\
Special case: To each $\delta\in \R^+$ exists $m_{\delta}$ so that the distance between adjacent zeros of $\gamma$ in $[m_{\delta}, \infty)$ is $\leq \delta \ \ (s_+=s_-)$.\\
Examples: $\gamma=\sin(t^2)$, or $(\sin(t^2))/(1+|t|).$\\

\textit{Proof of Proposition 5.4}: One can assume $\norm{u + \gamma v}_{\infty}=:a <\infty.$\\
We first show $\norm{v}_{\infty}<\infty$: Else there exist e.g. $t_m\to \infty $ with $\norm{v(t_m)}\to \infty, m\to\infty $ (the proof for $t_m \to -\infty$ is the same ). With the $\rho_0$ of (5.10) choose $\epsilon \in (0,1)$ with $\epsilon\norm{\gamma}_{\infty} <1, $ then $\delta \in (0, \rho_0)$ with
$\norm{u(s)-u(t)}\leq \epsilon$ and $\norm{v(s)-v(t)} \cdot \norm{\gamma}_{\infty} <\epsilon$  if $s,t \in J$ with $|s-t|\leq \delta$.
Choose $n\in \N$ with $\norm{v(t_m)}>\frac{2a+1}{\rho_0} + 3 $ if $m\geq n; $ with (5.10) to these $n, \delta$ there are $m_{\delta} \geq n $ and $s_+, s_- \in [t_{m_{\delta}}, t_{m_{\delta}} + \delta] $ with $|\gamma_+|\geq \rho_0$, $|\gamma_+  + \gamma_-| \leq \delta $, where $f_{\pm}:=f(s_{\pm})$. Then
\begin{align}
&\norm{\gamma_+v_+ - \gamma_- v_- } \leq \norm{u_+ +\gamma_+v_+} + \norm{u_--u_+} + \norm{u_-+ \gamma_-v_-} \\ &\leq
a+\epsilon+ a \leq 2a+1;\nonumber
\end{align} but
\begin{align*}
&\norm{\gamma_+v_+- \gamma_-v_-} = \norm{\gamma_+ ( v_++ v_-)- \gamma_++\gamma_-) v_-} \geq \norm{\gamma_+(v_++v_-)} - |\gamma_++\gamma_-| \norm{v_-} 		\\
&	\geq	\rho_0\norm{2v_-+(v_+-v_-)}-\delta\norm{v_-}\geq \rho_0(2\norm{v_-}-\norm{v_+-v_-})- \delta \norm{v_-} \geq (2\rho_0 - \delta)\norm{v_-}		\nonumber \\
&-\rho_0\epsilon \geq \rho_0(\norm{v_-}-\epsilon\geq\rho_0(\norm{v(t_{m_{\delta}}}-\norm{v(t_{m_{\delta}})-v_-}-\epsilon) \geq \rho_0(\norm{v(t_{m_{\delta}}}-2\epsilon) \geq			\nonumber \\
&\rho_0(\frac{2a+1}{\rho_0} + 3 - 2)=2a+1+\rho_0 >2a+1		.	\nonumber
\end{align*}
By the above we have now $a:=\norm{u+\gamma v}_{\infty}<\infty, b:=\norm{v}_{\infty}<\infty.$\\
If (5.9) is false, there exist $\epsilon_0>0$ and $t_m \in J, $ e.g. $t_m \to \infty$, with
$\norm{u(t_m)} \geq a + 3 \epsilon_0, m\in \N; $ to $\epsilon_0$ exist $\delta_0 \in (0,\epsilon_0)$ with $\norm{u(s)-u(t)}\leq \epsilon_0,
\norm{\gamma}_{\infty} \norm{v(s)-v(t)} \leq \epsilon_0 $ if
$s,t \in J, |s-t|\leq \delta_0.$ With $O_2$, to $(t_m)$ and $\delta:=\delta_0/(2+b)$ there exist $m_{\delta} $ and
$ s_+, s_- \in [t_{m_{\delta}},t_{m_{\delta}} + \epsilon]$ with $|\gamma(s_+)+\gamma (s_-)|\leq \delta .$ With $\delta\leq \delta_0$ one gets
\begin{align*}
&2a+ 6 \epsilon_0  \leq 2 \norm{u(t_{m_{\delta}})} \leq 2 ( \norm{u(s_+)} +\epsilon_0)=\norm{u(s_+)+u(s_+)}  +2 \epsilon_0 \leq		\\ \nonumber
&  \norm{u(s_+)+u(s_-)} + \epsilon_0 + 2\epsilon_0=\norm{u(s_+)+\gamma(s_+)v(s_+)+u(s_-)-\gamma(s_+)v(s_+)} + 3\epsilon_0 \leq	\\ \nonumber
& a + \norm{u(s_-)+\gamma(s_-) v(s_-)} + \norm{\gamma(s_+)v(s_+)+\gamma(s_-)v(s_-)} + 3\epsilon_0	\leq \\ \nonumber
& 2a+\norm{(\gamma(s_+) +\gamma (s_-))v(s_+)} + \norm{\gamma(s_-)(v(s_-)-v(s_+))} + 3\epsilon_0 \leq			\\ \nonumber
& 2a + \delta b + \epsilon_0 +3 \epsilon_0 \leq 2 a + \delta_0 + 4\epsilon_0 \leq 2 a + 5\epsilon_0.			
\end{align*}
Proof for $\gamma \in MC_0, v\in BUC:$ By the assumptions, one can apply Lemma 3.11, so $M_h(\gamma v) \in C_0$. With
$\norm{M_hu + M_h(\gamma v)}\leq a $ one gets $\limsup_{t\to\infty} |M_hu|\leq a;$
since $u\in UC, M_hu\to u$ uniformly on $J$, so $\limsup_{t\to\infty} |u| \leq a. $ If $J=\R$, similarly $\limsup_{t\to-\infty } |u| \leq a, $ therefore $\limsup_{|t|\to \infty} |u| \leq a$. --
\begin{corollary}
If $\gamma$ is as in Proposition 5.4, $U\subset PS_+(J,X)\cap UC(J,X)$ (see (3.6), $V\subset UC(J,X)$), then even
\begin{equation}
\norm{u}_{\infty} \leq \norm{u+\gamma V}_{\infty}, \  u \in U , v\in V.
\end{equation}
\begin{proof}
The definition of $PS_+$ gives $\norm{u}_{\infty} = \sup_J |u| = \limsup_{t\to \infty} |u| \leq \limsup_{|t|\to\infty } |u|, \leq \norm{u+\gamma v} $ by Proposition 5.4. -- \end{proof}
\end{corollary}
\begin{proposition}
If the Porada inequality (5.9) holds for $\gamma  \ ( \in \F^J), $ then, with this $\gamma, UC(J,X)\cap(\gamma UC(J,X)) \subset C_0(J,X)$.
If (5.9) holds for $\gamma $ with $u,v$ only $\in BUC(J,X)$, then, with this $\gamma$, $BUC(J,X)\cap(\gamma BUC(J,X)) \subset C_0(J,X)$.
\begin{proof}
If $u=\gamma v$ with $u,v \in UC$ resp. $BUC$, then $\norm{u-\gamma v}_{\infty} =0,$ so (5.9) gives $\limsup_{|t|\to\infty}|u|=0, u\in C_0$ follows with $u\in UC.$
-- \end{proof}
\end{proposition}
A converse of Proposition 5.8 is in general false: For $\gamma = 2 + g_{1,2}$ one has $UC\cap (\gamma UC) \subset C_0$ by Lemma 5.1, with Lemma 3.10 and (5.4) for $g_{1,2}$, but (5.9) is false ($u=2, v=-1$).
\begin{theorem}
If the Porada inequality $\limsup_{|t|\to\infty} |u|\leq \norm{u+\gamma v}_{\infty}$ holds for $u\in U, v\in V, U,V$ fixed $\subset X^J, \inf_J |\gamma | >0$, then $U+\gamma V$ is uniformly closed if $U$ and $V$ are uniformly closed, $U-U\subset U$ and $V-V\subset V$,
and if \emph{one} of the following 3 assumptions holds (with $PS_+$ of (3.6)):
\begin{enumerate}
\item[(a)] $U\subset PS_+(J,X)$.
\item[(b)] $V\subset PS_+(J,X),  \ \norm{\gamma}_{\infty}<\infty.$
\item[(c)] $ U + C_0(J,X) \subset U$ and $ V+C_0(J,X) \subset V, \ U,V \subset C(J,X),  \ \gamma \in C(J,\F)$,\\ $J$ closed.
\end{enumerate}
\end{theorem}
\begin{examples} $U,V \subset UC(J,X)$ and $\in \{X_c, C_0, X_c+C_0, AP, AAP, AA, AAA, $ \\ $LAP, A $ linear uniformly closed $\subset REC, A$ linear uniformly closed $\subset PS_+, BUC, UC\}$  ($12^2$ different pairs)) and $\gamma = \phi g_{\omega, r}$ as in Examples 5.5.
\end{examples}
This follows with Proposition 5.4 and Corollary 5.7; $BAA\cap UC= AA\cap UC; $ for $AA, LAP, REC$ see after (3.6).

\textit{Proof of Theorem 5.9, case (a):} Assume $(w_n)_{n\in\N}$ is a Cauchy sequence from $U+\gamma V$
with respect to $\norm{\ }_{\infty}, w_n=u_n+\gamma v_n, u_n\in U, v_n\in V; $
the $u_n, v_n $ are in general not unique; though $U$ and/or $V$ may contain unbounded functions, ``$(w_n)$ Cauchy'' is still defined:
If $\norm{w_n-w_m}_{\infty}= \norm{ (u_n-u_m)  + \gamma(v_n-v_m)}_{\infty}\leq \epsilon$ if $n,m\geq n_{\epsilon}$, with the assumed Porada inequality one gets $\limsup_{|t|\to \infty} |u_n-u_m| \leq \epsilon$ if $n,m \geq n_{\epsilon}.$
With (a), $u_n-u_m \in PS_+$, so $\norm{u_n-u_m}_{\infty} =\limsup_{t\to\infty}|u_n-u_m| \leq \epsilon$ by the proof of Corollary 5.7.
Since $U$ is uniformly closed, there exists $u\in U$ with $ \norm{u_n-u}_{\infty}\to 0$. Then $(w_n-u_n) = (\gamma v_n)$ is $\norm{ \ }_{\infty}$-Cauchy, with $\inf_J |\gamma| >0$ $(v_n)$ is Cauchy, there is $v\in V$ with $\norm{v_n-v}_{\infty} \to 0. \  w_n\to w:=u+\gamma v:$
The above implies $w_n \to w $ pointwise on $J, \norm{w_n-w_m}_{\infty} \leq \epsilon$ if $n,m \geq n_{\epsilon}$ and $n\to\infty$ gives
$\norm{w-w_m}_{\infty} \leq \epsilon,\  m\geq n_{\epsilon}$.\\
Case (b): $\norm{w_n-w_m} \leq \epsilon$ for $n,m\geq n_{\epsilon} $ and Porada give $\limsup_{t\to\infty} |u_n-u_m| \leq \epsilon$,
then $\limsup_{t\to\infty} |\gamma(v_n-v_m)| \leq \limsup_{t\to \infty} (|u_n-u_m|+|w_n-w_m|) \leq \limsup_{t\to\infty} |u_n-u_m| + \epsilon \leq 2 \epsilon;$ with $ \delta_0:=\inf_J |\gamma|$ one gets $\limsup_{t\to\infty} |v_n-v_m| \leq 2\epsilon/ \delta_0$, then $\norm{v_n-v_m} \leq 2\epsilon/\delta_0$ since $v_n-v_m\in PS_+$.
So $\norm{v_n-v} \to 0$ for some $v\in V$. With $\norm{\gamma}_{\infty}<\infty$ also $\norm{\gamma v_n-\gamma v}_{\infty} \to 0.$ $(w_n)$ Cauchy gives now $(u_n)$ Cauchy, then an $u\in U$ with $\norm{u_n-u}_{\infty} \to 0, $ so $\norm{w_n-w} \to 0$ with $w=u+\gamma v \in U+ \gamma V.$ -- \\

For case (c) we need
\begin{lemma}
To $u\in C(J,X)$ with $\limsup_{|t|\to\infty } |u| =: a <\infty$ and $J$ closed exists $\phi \in C_0(J,X)$ with $\norm{u+\phi}_{\infty}\leq a.$
\begin{proof}
$J\ne \R, a:=\limsup_{t\to\infty } |u|:$ With $v(t):=u(t)\cap a:=(\min \{\norm{u(t)},a\})$ $u(t)/\norm{u(t)}$ if
$u(t)\ne 0$, else $:=0,$ the $v\in C(J,X)$ with
$\norm{x\cap a - y\cap a } \leq 2 \norm{x-y}$ (\cite[p.327 (12)]{lit20}), the definition of $\limsup$ and $v$ give $\phi:=v-u \in C_0(J,X), |u+\phi| =|v| \leq a$ on $J$.\\
Case $J=\R$:\\ Applying the above to $w(t):=u(-t), t\geq 1$, one gets $\phi_- \in C_0((-\infty,-1]),X)$ with
$|u+\phi_-| \leq \limsup_{t\to - \infty } |u|\leq a$ on $(-\infty,1];$ the above for $u\big| [1,\infty) $ gives $\phi_+\in C_0([1,\infty),X)$ with $|u+\phi_+| \leq \limsup_{t\to\infty } |u| \leq a.$
Interpolating linearly between $-1$ and $1$ one gets $\Phi \in C_0(\R,X)$ with $|u+\Phi| \leq a$ on $\R$, also on $[-1,1]$, since $|u+\Phi| \leq a$ in $\pm 1$.
-- \end{proof}
\end{lemma}
\textit{Proof of Theorem 5.9, case (c)}: Assume $(w_n)$ a Cauchy sequence from $U+\gamma V, w_n=u_n+\gamma v_n, u_n \in U, v_n\in V, $ the $u_n, v_n$ in general not unique.
To given $\epsilon_n\in \R^+ $ with $\sum_1^{\infty} \epsilon_n <\infty$ one gets recursively a subsequence which we denote again by $(w_n)$ with $\norm{w_{n+1} - w_n} < \epsilon_n, n\in \N$.
With the assumed Porada inequality one has $\limsup_{|t|\to\infty } |u_{n+1}-u_n|\leq a_n:=\norm{w_{n+1}-w_n} < \epsilon_n, n\in \N$.\\
Define $z_1:=u_1\in U, \phi_1 :=0 \in C_0(J,X).$
With Lemma 5.11 there exists $\phi_2 \in C_0(J,X)$ with $\norm{z_2-z_1} = \norm{u_2-u_1+\phi_2} <\epsilon_1,
z_2:=u_2 + \phi_2 \in U$.\\
If one has already $z_1,\dots, z_{n+1}\in U$ with $z_k=u_k + \phi_k, 1 \leq k \leq n+1, \norm{z_{k+1}-z_k} < \epsilon_k, 1\leq k\leq n, $ with Lemma 5.11 there is $\psi_{n+2}\in C_0(J,X)$ with
$\norm{u_{n+2}-u_{n+1}+\psi_{n+2}}$ $< \epsilon_{n+1}, $ so $\norm{z_{n+2}-z_{n+1}}< \epsilon_{n+1}$ with $z_{n+2}:=u_{n+2} + \phi_{n+2} \in U, $ $\phi_{n+2}:=\phi_{n+1} +\psi_{n+2} \in C_0(J,X)$.\\
So by ``definition through recursion'' there exists a sequence $(\phi_n)$ from $C_0(J,X)$ with $\norm{z_{n+1}-z_n} < \epsilon_n, z_n:=u_n + \phi_n$, all $n\in \N$.
With $\sum_1^{\infty} \epsilon_n <\infty$ the $(z_n)$ is $\norm{ \ }_{\infty}$-Cauchy in $U$, so by assumption there is $u\in U$ with $\norm{z_n-u}_{\infty} \to 0.$\\
With $(w_n)$ Cauchy one gets $(w_n-z_n)=(w_n-(u_n+\phi_n))=(\gamma v_n - \phi_n)=\gamma( v_n -\frac{1}{\gamma}\phi_n)$ is Cauchy with $\frac{1}{\gamma}\phi_n \in C_0(J,X)$ by the assumptions on $\gamma$; then $v_n-\frac{1}{\gamma} \phi_n \in  V$, there is $v\in V$ with $\norm{v_n -\frac{1}{\gamma} \phi_n - v}_{\infty} \to 0$, as in the proof of (a) one gets
$\norm{w_n-(u+\gamma v)}_{\infty} \to 0, u+\gamma v \in U + \gamma V, \norm{\gamma}_{\infty} <\infty $ is not needed. --\\

\textbf{4.} $U+\gamma V\subset M(U+\gamma V)$\\
$U+\gamma V \subset M(U+\gamma V)$ is certainly true if $U\subset MU $ and $\gamma V\subset M(\gamma V)$.
For this see § 3, especially Propositions 3.1, 3.2, 3.12, 3.16 and the examples there, and Examples 3.18, 3.19, 3.20.
So for example one has \begin{align}
&U+\gamma V \subset M(U+\gamma V) \text{ if } \gamma=g_{\omega,r}, \omega \in \R, 1\leq r \in \R, \text{ and } \\
&U\in \{X_c,C_0, X_c+C_0, AP, AAP, AA, LAP, REC, PS_+, BUC, UC \}, \nonumber \\
&V\in \{C_0,AAP:=C_0 + AP, BUC, UC\}.\nonumber
\end{align}
(5.14) also holds if $C_0\subset U\subset MU, $ and only $V\subset BUC$, e.g. $V\in \{X_c, X_c+C_0, AP, AA\cap UC, LAP\cap BUC \}$ (though e.g. $gAP \not\subset M(gAP))$, so for ``asymptotic'' extensions $C_0 + gX_c, C_0+gAP=g(AAP), C_0 + g(AA\cap UC)$ (see \cite[Prop.3.4, Prop. 2.2(i), (3.3)]{lit9}).\\

\textbf{5.} $(\Delta)$ for $U+\gamma V$ is treated in § 6.
\abstand
\section{The property ($\Delta$) for vector sums}
Here we will discuss mainly $U+\gamma V$ with
\begin{equation}
U,V \in \{X_c,C_0,X_c+C_0,AP, AAP, BUC, UC\}.
\end{equation}
In the case $\gamma \equiv 1$ a first general result has been obtained in \cite[Proposition 5.1 p. 432]{lit9}:
\begin{theorem}
If $U$ and $V\subset \Lloc(J,X)$ are positive invariant additive groups with ($\Delta$), $J$ closed, any $X$, if $U\cap V =\{0\}, V\subset C(J,X)$ and there is $m_0 \in \N$ with $U\subset M^{m_0} Av_0(J,X)$, and finally $U$ and $V$ are invariant if $J=\R$, then $U+V$ also satisfies ($\Delta$).\\
\textnormal{Here for $J\ne \R$ resp. $J=\R$ (see \cite[p. 1007-1008]{lit8}, \cite[(3.4),(3.5) p. 42]{lit10})}
\begin{align}
Av_0(J,X):=&\{ f\in \Lloc(J,X): \text{ to } f \text{ exists } r_0 \in J \text{ with } \\
		&\frac{1}{T}\int_0^T f(r_0 + s) ds \to 0 \text{ as } T \to \infty\}, \nonumber\\
Av_0(\R,X):=&\{ f\in \Lloc(\R,X): \text{ to } f \text{ exists } \delta >0 \text{ with } \\
		&\frac{1}{2T}\int_{-T}^T f(r + s) ds \to 0 \text{ as } T \to \infty,\text{ uniformly in } r\in [0,\delta]\} \nonumber
\end{align}
\end{theorem}
Examples where Theorem 6.1 gives ($\Delta$) for $U+V$ are (see \cite[Corollaries 5.3, 5.8, Examples 5.5, 5.6, 5.9]{lit9})\\
$U+X_c $ with $U$ as in Theorem 6.1 and $U\cap X_c=\{0\}$, Banach space valued pseudo almost periodic functions $PAP$ and generalized pseudo almost periodic functions $GPAP$ (for range space $X=\C$ see \cite[p. 56, 66]{lit29}), asymptotic or pseudo almost automorphic functions $AAA:=C_0+AA$, $PAA:=Av_n + AA$, Eberlein weakly almost periodic functions $EAP$  ([13, p. 14/15]), analog extensions of Levitan almost periodic functions $LAP$ and of recurrent functions $REC$.
\begin{examples}
Theorem 6.1 gives $(\Delta)$ for $U+g_{\omega, r} V$ with $U\in \{X_c,AP,AA\}, V \in \{AAP, BUC\}$ if $J$ is closed, $0\ne \omega \in \R$ and $r=2$ for $V=AAP, 1<r \leq 2 $ for $V=BUC$.
Included are $AAP+gAP=AAP+gAAP, =AP+gAAP$.
\end{examples}
\textit{Proof}: $+$ is direct with Lemma 3.13 (a); $g_{\omega, r} BUC$ is positive invariant by § 5, 2(e), for ($\Delta$) see (4.15);($\Delta$) for $g_{1,2} AAP$ has been shown by B. Bolis and the author (2002, unpublished); for ($\Delta$) for $AA$ see \cite[Prop. 3.5 (ii)]{lit9}. --\\

$U+\gamma V$  with $U$ or $V=C_0$ is with $g_{\omega,r} C_0=C_0$ already included in the examples immediately after (6.3) resp. by § 4.\\

Theorem 6.1 is \emph{not} applicable already for $V=gX$ or $gUC$ (not positive invariant) or $gAP$ (has not ($\Delta$) by Example 3.14)) or $U,V \in \{AAP,BUC,UC\} $ (the sum $U+gV$ is not direct).\\

For cases not covered above see the matrix at the end of this §, e.g. $X+gAP, AP+gAP, $ and the table in § 8.\\

We turn now to the case $BUC + \gamma BUC$:\\
Here the following generalization of the Bohl-Bohr-Kadets criterion (see e.g. \cite[p. 674, 677 ($L_U$)]{lit7}) introduced by Loomis \cite[Theorem 3 p. 365]{lit25} turns out to be helpful:
\begin{definition}
We say that an $A \subset \Lloc(J,X)$ \emph{satisfies the Loomis condition } $(L_b)$ if (6.4) holds:
\begin{equation}
	f \in L^{\infty}(J,X) \text{ and all } \Delta_h f \in A, h\in \R^+ \text{, imply } f\in A.
\end{equation}
\end{definition}
By Proposition 3.12 of \cite[p. 682]{lit7}, if a linear positive invariant $A$ satisfies ($\Delta$) and $A\subset MA$, then $(L_b)$ is equivalent with the classical Bohl-Bohr criterion ($P_b$)
\begin{equation}
f \in A \text{ and indefinite integral Pf bounded implies } Pf\in A.
\end{equation}
\begin{proposition}
(a) If $A$ is linear $\subset BC(J,X), A\subset MA$ and $A$ satisfies $(L_b)$, then $A$ satisfies $(\Delta)$.\\
(b) If $BUC(J,X)\subset A\subset BC(J,X), A\subset MA$, $A$ is linear, positive invariant and satisfies ($\Delta$), then $A$ satisfies ($L_b$).
\begin{proof} (a): If $f\in \Lloc$ with all $\Delta_h f \in A \subset BC$, then $F:=f-M_h f \in BC$ since $BC$ satisfies ($\Delta$): By ($\Delta$) for $L^{\infty}$ of Proposition 1.1 of \cite[p. 677]{lit7}
the $F$ is a.e. bounded, by Lemma 2.1 $f$ and so $F$ are continuous,
so $F\in BC$. Since $\Delta_k M_h f=M_h\Delta_k f \in A$ with $A\subset MA,\ \Delta_k F \in A$: with ($L_b$) one gets $f-M_hf = F \in A, h\in \R^+.$\\
(b) If $f \in A$ with $Pf$ bounded, $Pf \in BUC \subset A$, so $(L_b)$ for $A$ follows with the above mentioned Proposition 3.12 of \cite{lit7}.
-- \end{proof}
\end{proposition}
Proposition 6.4 (b) becomes false if only $A \subset BUC: A=AP(\R, c_0(\N,\F))$ is linear invariant with
($\Delta$) and $A\subset MA$, but $(L_b)$ is false by Amerio \cite[p. 53/54]{lit1}
($A\subset MA, (L_b) \Rightarrow (P_b).$)\\

With Proposition 6.4 (a), to get ($\Delta$) for $BUC+\gamma BUC$, we prove ($L_b$);
since besides $f\in \Lloc$ with all differences $\Delta_hf \in A$, one has the additional information ``$f$ bounded'', which seems to make matters a bit easier:\\

Assumptions and notation in all of the following, till Theorem 6.5, will be

$X$ a Banach space with scalar field $\F, (\Delta_h f)(t):=f(t+h)-f(t)$
\begin{align}
&J=\R, BUC:=BUC(\R,X), A:= BUC + \gamma BUC, BC:=BC(\R, X), \\ \nonumber  &C_0:=C_0(\R,X)\\
&\gamma \in BC(\R, \F), \inf_{ } {}_{\R} |\gamma|>0, \gamma BUC \text{ invariant }, \gamma BUC \subset M(\gamma BUC), \\
&\text{the Porada inequality (5.9) holds for } \gamma \text{ and only } u,v \in BUC.\nonumber \\
&F:\R\to X \text{ Bochner-Lebesgue measurable with } \norm{F}_{\infty}:= \\
&\sup \{\norm{F(t)}: t\in \R\} <\infty \text{ and } \Delta_h F \in A \text{ for all } h \in \R^+.\nonumber
\end{align}
Substituting $t-h$ for $t$ in $\Phi:= \Delta_h F \in A$, one gets $\Delta_{-h}F=-\Phi_{-h}, \in A$,
since the invariance (2.7) of $\gamma BUC$ by (6.7) implies that of $A$,
$\Delta_{-h}F \in A$, i.e.
\begin{align}
&\Delta_h F \in A \text{ for all } h\in \R
\end{align}
(6.6), (6.7) and (6.8) give with Lemma 2.1
\begin{align}
& F\subset BC.
\end{align}
With the $\norm{ \ }_{\infty} $ as norm, $C_0, BUC, \gamma BUC, A, BC$ are all Banach spaces, $C_0, BUC,$ $\gamma BUC, A$ are closed linear subspaces of $BC$, with $C_0 \subset BUC$, $C_0\subset \gamma BUC$ by the assumptions for $\gamma$ in (6.7), $C_0 \subset A \subset BC$.\\
The factor spaces
\begin{equation}
U:=BUC/C_0, \ V:=(\gamma BUC)/C_0, \ Y:=A/C_0, \ Z:=BC/C_0
\end{equation}
with e.g. $U=\{[x]: x\in BUC\}, [x]:=x+C_0, =u+C_0=[u]$ for all $u\in [x]$,
are all well defined, with $[u]+[v]:=[u+v], r[x]:=[rx]$ and the norm $\norm{ [x]}:= \inf \{\norm{u}_{\infty}: u\in [x]\}$ they all become Banachspaces, with
$$U\subset Y, V\subset Y, U+V\subset Y \subset Z,$$
$U$ and $V$ resp. $Y$ are closed linear subspaces of $Y$ resp. $Z$ (see e.g. Day \cite[Ch. II § 1 Lemma 1 p. 29]{lit16}).

The assumed Porada inequality of (6.7) and Proposition 5.8 give $BUC\cap (\gamma BUC) \subset C_0, $ so
\begin{equation}
U\cap V=\{0\}, Y=U + V = U \oplus V \text{ is a direct sum}
\end{equation}
The translation $BC \ni f \to f_t \in BC, (f_t)(s):=f(s+t), s,t \in \R$ is transferred in the canonical way to $Z$:
\begin{align}
&\text{For } t\in \R, T(t):Z\to Z \text{ is defined by } T(t)z:=[f_t], f\in z\in Z, \end{align}
 independent of $f\in z.$

Then $T(t): Z\to Z$ is linear, with $\norm{f_t}_{\infty}=\norm{f}_{\infty}, t \in \R, $
\begin{equation}
\norm{T(t)z}=\norm{z}, \  z\in Z, t\in \R;
\end{equation}
$T(t)$ is continuous on $Z$, so
\begin{equation}
T(t) \in L(Z,Z):=\{S:Z\to Z: S \text{ linear continuous} \}, \norm{T(t)}=1, t\in \R.
\end{equation}
The definition gives also
\begin{align}
&T(s)T(t)=T(s+t),  \ s,t \in \R, T(0)z=z, z\in Z, \\
&T: \R \to L(Z,Z) \text{ is a one-parameter group}\nonumber\\
&T|U: U\to U \text{ is even a } C_0 \text{ group},
\end{align}
i.e. (6.15), (6.16) with $U$ instead of $Z$ hold and $\norm{T(t)u-u}\to 0$ as $0<t\to 0, u \in U:$
this follows with the uniform continuity of $f\in u, \norm{f_t-f}_{\infty} \to 0, t \to 0.$ Even
\begin{equation}
T( \cdot ) u \in BUC(\R,U) \text{ for any } u \in U
\end{equation}
with (6.14), (6.15), (6.16) (this holds even for $C_0$-semigroups e.g. by \cite[p. 4 Cor. 2.3]{lit26}).\\
The invariance of $\gamma BUC$ of (6.7) gives analogously with (6.15)
\begin{equation}
T(t)V\subset V,  \ t\in \R, T(t)|V \in L(V,V)
\end{equation}
Now with (6.9), (6.13) and $z:=[F]\in Z$ one has $T(s)z-z=:w(s)\in Y=U+V$ of (6.12), so
\begin{equation}
T(s)z-z=:u(s)+u(v),  \ s\in \R, \text{ with unique } u(s)\in U, \  v(s)\in V.
\end{equation}
So with (6.16)
\begin{align*}
&T(s+h)z-z=T(s+h)z-T(s)z+T(s)z-z = \\&T(s)(T(h)z-z)+T(s)z-z=
T(s)(u(h)+v(h)) + u(s) + v(s),
\end{align*}
so
\begin{equation}
u(s+h)-u(s)-T(s)u(h)=-v(s+h)+v(s)+T(s)v(h), \ s,h\in \R,
\end{equation}
Since $T(s)U\subset U$ and $T(s)V\subset V$ by (6.17), (6.19) and the $U,V$ are linear with $U\cap V=\{0\}$ by (6.12), one gets
\begin{equation}
u(s+h)-u(s)=T(s)u(h), \  s,h\in \R
\end{equation}
with $u:\R \to U$ defined by (6.20) and $T( \cdot) u(h)\in C(\R,U)$ for any fixed $h\in \R$ by (6.18); similarly (but $T(\cdot)v(h)$ is in general not continuous)
\begin{equation}
v(s+h)-v(s)=T(s)v(h), \  s, h \in \R
\end{equation}
If $u$ were $\in \Lloc $ or even only Bochner-Lebesgue measurable one would get $u\in C(\R, U)$ with Lemma 2.1, but no such property can be deduced from (6.20).

In fact, for any additive $\beta: \R \to \F$ (and there are non-measurable ones, e.g. the linear functional defined by the hyperplane in Aufgabe 117 of \cite[p. 189]{lit20}) and any constant $x_0 \in X_c \subset BUC$ the
$u(s):=\beta(s)x_0$ defines a solution of the functional equation (6.22).

To get regularity properties for $u$, we need first boundedness of the $u$ of (6.20): For fixed $s\in \R$ one has with (6.20) and $z=[F]$
$$[F_s-F]=[F_s]-[F]=T(s)z-z=u(s)+v(s)=[p]+[q]=[p+q]$$
with $p\in u(s), q\in v(s).$ So there is $r\in C_0$ with $p+q=F_s-F+r $ or
$\norm{p+(q-r)}_{\infty} = \norm{F_s-F}_{\infty} \leq 2\norm{F}<\infty$ by (6.8), with $p\in BUC, q-r\in \gamma BUC + C_0=\gamma BUC$
with (6.7). With the assumed Porada inequality of (6.7) one gets
$$\limsup_{|t|\to \infty} |p| \leq 2\norm{F}_{\infty} , [p]=u(s);$$
With Lemma 5.11 there is $\phi \in C_0 $ with $\norm{p+\phi}_{\infty}\leq 2\norm{F}_{\infty}$, and then, $s$ being arbitrary $\in \R$,
\begin{equation}
	\norm{u(s)} \leq 2 \norm{F}_{\infty}, s \in \R.
\end{equation}
 If now $y\in \text{dual }U'$, with $\Phi(s):=y(u(s)), s\in \R$, with (6.22) one gets
$\Phi: \R \to F$ with
\begin{equation}
\Phi(s+h) - \Phi(s)=\Psi(s,h):=y(T(s)u(h)), s\in \R, h\in \R,
\end{equation}
with $\Psi(\cdot, h)\in BC(\R,\F)$ for each fixed $h\in \R$ by (6.18).
Here one can now apply a result of De Bruijn \cite{lit17}, formulated and extended in § 7; Theorem 7.2 there gives:
\begin{align}
&\text{To } \Phi \text{ of  (6.25) exists } G\in C(\R,\F), H \text{ additiv }: \R \to F \text{ with}\\
&\Phi=G+H \text{ on }\R; \nonumber
\end{align}
here $H$ additive means $H(s+t)=H(s)+H(t), s,t \in \R$.\\
In general, depending on $\Phi$, the $H$ need not even be measurable (see the example after (6.23)). Here however (6.24) and (6.26) imply the boundedness of $H$ on $[-1,1]$; Lemma 7.1 can be applied, yielding $H(s)=a_0 s, s\in \R, $ with some $a_0\in \F$, so
\begin{equation}
\Phi \in C(\R,\F).
\end{equation}
Since $y\in U'$ was arbitrary, the $u$ of (6.20), (6.22) is weakly continuous on $\R$.

With $\Q:=$ rational field $\subset \R$, the countable $u(\Q)$ is weakly dense in $u(\R)$,
with $\Q$ independent of the $y\in U'$, i.e. $u(\R)$ is weakly separable; this implies $u(\R)$ is (norm-)separable in $U$ (see e.g. \cite[p. 337 Lemma 1]{lit20}); furthermore $u$ is weakly (Lebesgue-)measurable since it is weakly continuous. Now Pettis's theorem (\cite[p. 7 Theorem 1.1.1]{lit2}) can be applied, $u$ is Bochner measurable; with (6.24) this gives
\begin{equation}
u\in \Lloc (\R, U).
\end{equation}
For the following we need
\begin{equation}
W:=\{\Phi \in U^{\R}: \text{ to } \Phi \text{ exists } z\in U \text{ with } \Phi(s)=T(s)z, s\in \R\}.
\end{equation}
$W$ is a linear invariant uniformly closed subspace of $BUC(\R,U)$ by (6.15), (6.16), (6.18), so $W$ satisfies ($\Delta$) by Proposition 4.1. $\Delta_h u \in W$ by (6.22) means $u\in \Delta W; $ Proposition 4.6 can be applied, so there exist $a,b \in U$ with $(u(0)=0 $ with (6.20), $T(0)a=a$ by (6.16))
\begin{equation}
u(s)=T(s)a - a + \int_0^s T(t)bdt,  \ s\in \R.
\end{equation}
$a=[\phi]$ with some $\phi\in BUC$, then $T(s)a-a=[\Delta_s \phi].$\\
Now to $\int_0^s T(t) b dt$ of (6.30): $b=[\psi]$ with $\psi \in BUC, H(\cdot) := T(\cdot)b \in W \subset BUC(\R,U)$.
So to $\epsilon>0$ exists $n\in \N$ with $\norm{\psi(t+h) -\psi(t)}\leq \epsilon$ and \\$\norm{H(t+h)-H(t)}\leq
\epsilon$ if $|h|\leq |s|/n, t\in\R$. For fixed $s\in\R$ define $s_j:=(js)/n, 1\leq j\leq n;$ then
$\norm{ \int_0^s T(t) b dt - \sum_1^nT(s_j) b\frac s n} \leq |s|\epsilon, \norm{ \int_0^s \psi(t+r)dr - \sum_1^n \psi(t+s_j)\frac s n } \leq |s| \epsilon$ for all $t\in \R$ (the first Bochner-integral for $U$-valued functions, the second for $X$-valued functions). Since $T(s_j)b=T(s_j)[\psi] =[\psi_{s_j}]=[\psi( \cdot + s_j)], $ one gets
\begin{equation}
\int_0^s T(r)bdr=[\Phi (s)], \text{ with  } \Phi(s)(t):= \int_0^s \psi(t+r) dr, \ t\in \R,
\end{equation}
with $\Phi \in UC(\R,BUC).$\\
With $\chi (t):=\int_0^t \psi(r)dr, \chi \in C(\R,X)$, one has $\Delta_s \chi(t)=\int_t^{t+s} \psi(r)dr=$ \\$\int_0^s \psi(t+r) dr , \ t \in \R, $ so
\begin{equation}
\Delta_s \chi=\Phi (s), \in BUC,  \ s\in \R.
\end{equation}
(6.30) and $T(s)a-a=[\Delta_s \phi]$ yield then
\begin{equation}
u(s)=|\Delta_s (\phi + \chi)|,\  s\in \R, \ \phi, \psi \in BUC, \chi=P\psi.
\end{equation}
Define
\begin{equation}
G:=F-\phi - \chi , \in C(\R,X),
\end{equation}
then with (6.20) and (6.33) one has, for $h\in \R$,
$$[\Delta_h G]=[\Delta_h F]-[\Delta_h(\phi + \chi)] = u(h) + v(h) - u(h)=v(h)\in V,$$
implying
\begin{equation}
\Delta_h G \in \gamma BUC, \ h \in \R
\end{equation}
Now $\gamma BUC \subset M(\gamma BUC), \gamma \in C(\R,\F)$ and $\inf_{\R} |\gamma| >0$ by the assumptions
(6.7), so Theorem 4.3 (b) can be applied, $\gamma BUC$ satisfies ($\Delta$).
With (6.35) there exists therefore to each $h\in \R^+$ a $\lambda (h)\in \gamma BUC$ with $G-M_hG=\lambda(h); $ (6.34) gives then
\begin{equation}
	F-M_hF-(\phi-M_h \phi) - (\chi - M_h \chi) = \lambda (h),\ h\in \R^+.
\end{equation}
Here $\phi \in BUC$ by (6.33), so $M_h\phi \in BUC$; (6.32) and ($\Delta$) for $BUC$ by Proposition 4.1 give $\chi-M_h \chi \in BUC, \ F\in BC$ by (6.10) implies $M_hF \in BUC$, so (6.36) gives
\begin{equation}
F\in BUC + \gamma BUC
\end{equation}
We have now shown
\begin{equation}
\text{If } J=\R \text{ and } \gamma \text{ satisfies (6.7), } BUC + \gamma BUC \text{ satisfies } (L_b).
\end{equation}
Now Proposition 6.4 (a) can be applied ($\gamma BUC \subset M(\gamma BUC)$ of (6.7) implies $A:=BUC + \gamma BUC \subset MA$):
\begin{equation}
\text{If } J=\R \text{ and } \gamma \text{ satisfies (6.7), } BUC + \gamma BUC \text{ satisfies } (\Delta).
\end{equation}\\

Extension of (6.39) to $J=[\alpha,\infty)$, with $\gamma \in BC(\R, \F)$ satisfying (6.7):
If $F\in \Lloc (J,X)$ with $\Delta_h F=u + (\gamma|J) v \in A:=BUC(J,X) + (\gamma | J) BUC(J,X),$ then
$F\in C(J,X)$ by Lemma 2.1. Define $G:=F $ on $J, G:=F(\alpha) $ on $(-\infty, \alpha),$ then
$G\in C(\R,X)$; with fixed $h \in \R^+ $ define further $U:=0$ on $(-\infty, \alpha-h]$, $U=u$ on $[\alpha,\infty), U$ linear on $[\alpha-1,\alpha]$, similar for $V$ with $v$;
then $U,V \in BUC(\R,X), U+\gamma V \in A(\R,X).$  $\phi:=\Delta_hG - (U+\gamma V)=0$ on $J=[\alpha,\infty)$ and
on $(-\infty, \alpha-h), $ so $\phi=0$ on $(-\infty, \alpha-h] \cup J; \ \  \phi \in C(\R,X)$ since $G,U,V, \gamma \in C(\R, \cdots)$, so $\phi \in C_0(\R,X) \subset BUC(\R,X)$, with $U^* := U+\phi$ one has $\Delta_h G= U^* + \gamma V \in A(\R,X). $ So (6.39) is applicable, yielding $G-M_h G \in A(\R, X)$ and then $F-M_hF = G|J - (M_h G)|J = (G-M_hG)|J \in (A(\R,X))|J=A(J,X)$:
\begin{align}
&\text{If } \alpha \in \R, J=[\alpha, \infty), \  \gamma| \R \text{ satisfies (6.7),} \\ &\text{then } BUC(J,X) + (\gamma| J) BUC(J,X) \text{ satisfies } (\Delta). \nonumber
\end{align}
Extension to $J=(\alpha,\infty), $ with $\gamma $ as in (6.40): If $F\in \Lloc( J,X)$ with $\Delta_h F \in A(J,X)\subset C(J,X)$, Lemma 2.1 gives $F\in C(J,X)$. If $h\in \R^+$ and $\Delta_h F(t)=F(t+h)-F(t)=u(t,h) + \gamma (t) v(t,h)$ with $u( \cdot, h), v(\cdot, h)\in BUC(J,X), u(\alpha + , h):=\lim_{\alpha < t \to \alpha} u(t,h) $ exists for each $h\in \R^+$, similarly $v(\alpha +, h)$, and $U(t,h):=u(t,h)$ for $t\in J, U(\alpha, h):=u(\alpha+, h)$ gives $U(\cdot , h) \in BUC([\alpha,\infty), X)$, similarly for $V(\cdot, h)$, with $F(t+h)-F(t)=U(t,h)+\gamma (t) V(t, h)$ for $t\in J, h\in \R^+; $ continuity of $F$ in $\alpha +h$ and $U,V$ in $\alpha$ gives
$F(\alpha + ) := \lim_{\alpha< t \to \alpha} F(t)$ exists $\in X$, independent of $h$, with $F^* (t):= F(t)$
if $t\in J, F^*(\alpha):=F(\alpha+)$ gives then $\Delta_hF^*(t)=U(t,h)+ \gamma (t) V(t,h)$ for $t\in [\alpha,\infty), h\in \R^+:$
One can apply (6.40) to $F^*$, getting $F^* - M_hF^* \in A([\alpha,\infty),X),$ restriction to $J$ gives $F-M_h F \in A(J,X)$.\\

Collecting all the above, we have shown
\begin{theorem}
If, for a Banach space $X, \gamma \in C(\R,\F_X)$ satisfies (6.7), then\\
$BUC(J,X) + (\gamma |J ) BUC(J,X)$ satisfies $(\Delta)$ for any $J$ as in § 2.
\end{theorem}
\begin{corollary}
If for $J$ and $X$ as in § 2 with $\F_X=\C$ the $\gamma $ satisfies
\begin{align}
&\gamma \subset BC(J, \F), \inf_J |\gamma| >0, \gamma \in MC_0 (J,\F), \text{ there exists } \epsilon_0 \in  \R^+ \text{ with  }
\\ \nonumber  &\gamma_a/ \gamma \in BUC(J,\F) \text{ if } 0 <a <\epsilon_0,
\end{align}
then $BUC(J,X)  + \gamma BUC(J,X)$ satisfies  ($\Delta$).

\begin{proof}
If $J=[\alpha,\infty), $ define $y:= \gamma $ on $[\alpha,\infty), y:=g=g_{1,2} $ on $(-\infty, \alpha-1], y $
continuous on \R and $\ne 0$ on $[\alpha-1 , \alpha] $ (e.g. piecewise linear).
Then $y \in BC(\R,X), \inf_{\R} |y| >0, y\in MC_0(\R, \C)$ with (6.41) resp. Lemma 3.10.
If $0<a<\epsilon_0, y_a/y \in BUC(J,\C)$ resp. $y_a / y =g(a) e^{2iat}$ on $(-\infty, \alpha - 1]$, so $\in BUC((-\infty , \alpha - 1], \C)$; since $y_a/y \in C(\R, \C), \phi:=y_a/ y \in BUC(\R, \C), 0<a<\epsilon_0;$
with $\inf_{\R} |\phi| \geq (\inf_{\R}|y|)/\norm{y}_{\infty} >0$ also $1/\phi \in BUC(\R,\C)$, then $1/\phi_{-a} \in BUC$,
i.e. $(y_{-a})/y \in BUC, $ so $y_a/y \in BUC(\R,\C)$ if $a\in (-\epsilon_0, \epsilon_0)$.
With the Example $BUC$ to Lemma 5.3, $yBUC(\R,X)$ is invariant. \\
Finally $y\in MC_0(\R, \C)$ with $(6.41)$ resp. Lemma 3.10, so by Proposition 5.4 the Porada inequality (5.9) holds for $u,v \in BUC(\R, X)$, and $yBUC(\R,X) \subset M(yBUC(\R,X)) $ by Proposition 3.12.\\
$y$ fulfills therefore all conditions of (6.7), one can apply Theorem 6.5, getting $(\Delta)$ for $BUC(J,X) + \gamma BUC(J,X)$ by restriction.

If  $J=(\alpha, \infty)$, by (6.41) $u:=\gamma_a /\gamma \in BUC(J,\C)$ if $0 <a< \epsilon_0$,
so $u(\alpha + ) := \lim_{\alpha<t \to \alpha} u(t)$ exists $\in \C, u(\alpha + )\ne 0$ with $0<  \inf_J |\gamma |, \sup_J|\gamma| <\infty$ of (6.41);
with the continuity of $\gamma $ at $\alpha + a$ then also $\gamma (\alpha + ):=\lim_{\alpha<t \to \alpha} \gamma (t)$ exists $\in \C$\textbackslash$\{0\}$. The extension $y$ of $\gamma$ to $[\alpha,\infty)$ by $y(\alpha):=\gamma (\alpha + )$ satisfies then (6.41) for $J=[\alpha,\infty)$. By the above case of closed $J$ now $A:=BUC([\alpha, \infty),X) + yBUC([\alpha,\infty),X)$ satisfies ($\Delta$).\\
If now $F\in L_{\text{loc}}^1 ((\alpha,\infty), X)$ has $\Delta_h F \in A_0:=BUC((\alpha,\infty), X) + \gamma BUC((\alpha, \infty), X)$ for $h\in \R^+$, by the proof of the case $J=(\alpha, \infty)$ for Theorem 6.5 the $F(\alpha+) $ exists and the corresponding extension $F^*$ of $F$ satisfies $\Delta_h F^* \in A, h>0$. $A$ has ($\Delta$), so $F^*-M_hF^* \in A $ for $h>0$, restriction to $(\alpha,\infty)$ gives $(\Delta)$ for $A_0$.
-- \end{proof}
\end{corollary}
\begin{remark}
There exists no real-valued $\gamma $ satisfying (6.41).
\end{remark}
\begin{examples}
For any $J$, complex $X, 0\ne \omega \in \R, 1\leq r \leq 2, BUC(J,X)+g_{\omega,r} BUC(J,X)$ satisfies ($\Delta$).
\begin{proof} $r=1:$ $g_{\omega,1} BUC=BUC;$ $1<r\leq 2$: $g_{\omega,r}$ satisfies (6.41) by Lemma 3.10, Proposition 3.12 and example $BUC$ in § 5, 2(e), Corollary 6.6 applies.
-- \end{proof}
\end{examples}
\begin{corollary}
If $\gamma $ satisfies (6.7) and $J$ is arbitrary, or only $\gamma | J$ is given satisfying (6.41), then $BUC(J,X) + \gamma BUC(J,X)$ satisfies $(P_b), (L_b).$
\begin{proof}
Proposition 6.4 (b) and Theorem 6.5 resp. Corollary 6.6.
-- \end{proof}
\end{corollary}
\begin{proposition}
If $U\subset UC(J,X)$ is linear and positive invariant, then \\
$\Delta(U+gX) = \Delta U, \ g(t)=e^{it^2}$.
\begin{proof}
$U\cap gX \subset BUC \cap gPS_+=\{0\}$ by (3.7). So if $f\in \Lloc$ with $\Delta_h f \in A:=
U+gX, h\in \R^+,$ then $\Delta_h f=u(h)  + ga(h)$ with unique $u(h) \in U, a(h)\in X, h\in \R^+$.
With $\Delta_{h+k}f = (\Delta_h f)_k + \Delta_k f$ one gets
$\Phi:=u(h+k)- u(h)_k -u(k)=g(g(k)e^{2ikt} a(h) + a(k)-a(h+k)), h, k \in \R^+$;
since the right side is bounded, $\Phi \in BUC,$ $(\cdots)$ on the right side $\in AP,$
again with (3.7) one gets $g \cdot (g(k)e^{2ikt} a(h) + a(k) - a(h+k))=0$ on $J$.
$d/dt$ gives $a(h)=0, h\in \R^+, $ and so $\Delta_h f \in U.$
-- \end{proof}
\end{proposition}
So if $U$ is as in Proposition 6.10
\begin{equation}
U + gX \text{ has } (\Delta), \text{ if } U \text{ has } (\Delta).
\end{equation}
\begin{examples} [to (6.42)] $U=\{0\}, X_c, C_0, X_c + C_0, AP, AAP, BUC, UC.$
\end{examples}
A refinement of the proof of Proposition 6.10 shows
\begin{remark}
Proposition 6.10 holds more generally for $U+\gamma X$ with $\gamma $ satisfying $O_1, \norm{\gamma}_{\infty}<\infty$
and for which exist $a\in \R^+, t_0 \in J$ so that $\gamma (t) \ne 0$ for $t>t_0$ and $\gamma_a /\gamma $ is not constant on $[t_0,\infty)$.\\
Special case: $\gamma \in MC_0, \norm{\gamma}_{\infty} <\infty, \liminf_{|t|\to \infty} |\gamma| >0;$
the $g_{\omega,r}$ with $0 \ne \omega \in \R, 1<r\in \R$ satisfy these conditions.
\end{remark}
Proposition 6.10 becomes however in general false for $U+g_{\omega,r}X$ if $r=1$ (also  if $\omega=0$). See also Example 4.4.
\begin{corollary}
If $0\ne \omega \in \R, 1<r\leq 2 , \ \ U, V \in \{BUC(J,X), UC(J,X)\}$ (4 pairs),
then for any $J$ and $X$ the $U+g_{\omega, r} V$ satisfy ($L_b$) and ($P_b$).
\begin{proof}
$U=V=BUC$: Corollary 6.9 with proof of Example 6.8.
In the other 3 cases, if $f\in L^{\infty} (J,X)$ and $\Delta_h f = u + g_{\omega, r} v \subset A:= U+g_{\omega, r} V$
for all $h\in \R^+$, the $f \in C(J,X)$ by Lemma 2.1 and $\Delta_h f \in BC, $ so $\norm{u+g_{\omega,r} v } _{\infty}< \infty$. By Example 5.5 the
Porada inequality (5.9) holds,
yielding $\limsup_{|t|\to \infty} |u| <\infty$; since $u\in UC$ one gets $u \in BUC$, then $v\in BUC$, $\Delta_h f \in BUC + g_{\omega, r} BUC.$
The proof of Example 6.8 gives (6.41), so Corollary 6.9 yields $f\in BUC + g_{\omega, r} BUC \subset A$, i.e. ($L_b$). ($P_b$) follows from ($L_b$), since $A\subset MA $ here
with (4.15).
-- \end{proof}
\end{corollary}
\abstand
\begin{center}
\textbf{$(\Delta)$ for some further vector sums $U+gV$}\\
\noindent $g(t)=e^{it^2}, \ J=\R, \ X$ complex Banach space
\end{center}
\begin{table*}[h]
	\centering
		\begin{tabular}{|c|c|c|c|c|c|c|}
			\hline \backslashbox{$U$}{$gV$} 			& $gX$ & $g(X_c+C_0)$ & $gAP$ & $gAAP$ & $gBUC$ & $gUC$	 \\
			\hline $X_c$ 		& $+$ & $+$ & $-$ & $+$ & $+$ & $+$ \\
			\hline $C_0+X_c$ & $+$ & $+$ & $+$ & $+$ & $+$ & $+$	 \\
			\hline $AP$ 			& $+$ & $+$ & $-$ & $+$ & $+$ & $+$	 \\
			\hline $AAP$ 		& $+$ & $+$ & $+$ & $+$ & $+$ & $+$	 \\
			\hline $BUC$ 		& $+$ & $+$ & $+$ & $+$ & $+$ & $+$ \\
			\hline $UC$ 			& $+$ & $+$ & $?$ & $?$ & $?$ & $?$	\\
			\hline
		\end{tabular}
	\label{tabelle1}
\end{table*}

\noindent\textbf{Comments} (j:k refers to row j, column k)
\begin{itemize}
\item[+] means the $U+gV$ here satisfies ($\Delta$)
\item[-] means the $U+gV$ here does not satisfy ($\Delta$)
\item[?] means the $(\Delta)$ for this $U+gV$ is still open
\end{itemize}
\textbf{Identities}: $gX=gX_c,  gC_0=C_0,   C_0 +g X=g(X_c + C_0),$\\
$C_0 + X_c + gAP=X_c + gAAP = X_c+C_0 + gAAP, C_0+gAP=C_0+gAAP=gAAP,$\\
$C_0+gBUC=gBUC, C_0+gUC=gUC, AP+gAAP=AAP+gAP=$\\
$AAP+gAAP,$\\
$AAP+gBUC=AP+gBUC, AAP+gUC=AP+gUC,$\\
$BUC+gAAP=BUC+gAP,UC+gAAP=UC+gAP.$\\

\noindent\textbf{Proofs}:
\begin{itemize}
\item[column 1:] (6.42), Proposition 4.1.
\item[column 6:] j:6 can be reduced to j:5, $1\leq j \leq 5$, using the (any $J,X, g(t)=e^{it^2})$ \textbf{Proposition 6.14}: $U$ linear positive-invariant $\subset BUC$ implies \\$\Delta(U+gUC)=\Delta(U+gBUC)$.
\item[row 1:   ] 1:4, 1:5 Theorem 6.1; 1:6 see column 6.\\
									 1:2=2:1; 1:3 contradiction with Prop. 3.1, Lemma 3.11, (3.7).
\item[row 2:   ]  2:2=2:1; 2:3=2:4=1:4; 2:5=1:5; 2:6=1:6 (see after (6.3)).
\item[row 3:   ] 3:2=4:1; 3:3 as 1:3; 3:4, 3:5 Theorem 6.1, Lemma 3.13 (a).
\item[row 4:   ] 4:2=4:1, 4:3=4:4=3:4; 4:5=3:5.
\item[row 5:   ] 5:2=5:1; 5:5 Ex. 6.8; 5:4=5:3: $\Delta_h F(t)=u(t,h)+g(t) v(t,h)$ uniquely, Prop. 7.3 gives $G=\phi + H$ for $G(t):=u(\beta, t-\beta), H$ additive; $M_nM_1(\Delta_h(F-\phi))  \to H(h)$ pointwise, so $H$ is linear, then $G-M_hG \in BUC; (\Delta)$ for $BUC +gAP$ follows with $\Delta gAP=gAP+X_c $ (B.Basit 2008).
\item[row 6:   ] 6:2=6:1, 6:3=6:4.
\end{itemize}
For 1:4, 3:4 ($\Delta$) for $gAAP$ is needed: See the proof of Ex. 6.2. \\
Since the $gX$ is not positive invariant and $gAP$ has not ($\Delta$), in these cases Theorem 6.1 is not applicable.
\abstand
\section{On functions with continuous differences}
\begin{lemma}
If $J=\R$ or $\R_+$, $X$ a Banach space, $f:J\to X$ satisfies
\begin{equation}
f(s+t)=f(s)+f(t) \text{ for all } s, t \in J,
\end{equation}
and $f$ is bounded on some non-empty open set $\subset J$, then there is $a \in X$ with $f(s)=sa, \ s\in J$.\\
\textnormal{For generalizations see \cite[e.g. Theorem 4.1]{lit3}.}
\begin{proof}
With the assumptions and translation one can assume that $|f|\leq \beta $ on $[1-\delta, 1]$ with some $\beta,  \delta \in \R^+$. If $f$ is not continuous
from the left in some $s_0 \in J$, with translation one can ssume $s_0=1;$
so there are $s_n \in [1-\delta,1)$ with $s_n \to 1 $ but $\norm{f(s_n) - f(1)} > \epsilon_0$ for some
$\epsilon_0 \in \R^+$. To $\epsilon_0$ and $\beta$ exists $m_0 \in \N$ with $m_0\epsilon_0 >2 \beta$; then to $m_0$ and $\delta$
exists $n_0 \in \N$ with $ 1 - s_{n_0} < \delta/m_0;$ then $1-m_0 (1-s_{n_0}) \in [1-\delta, 1]$; with
$a:=m_0(1-s_{n_0}) \in \R^+$ and $f(1-a) + f(a)=f(1)$ one gets
\begin{align*}&2\beta \geq \norm{f(1-a)-f(1)} = \norm{f(a)}=m_0\norm{f(1-s_{n_0})}=m_0 \norm{f(1)-f(s_{n_0})}\\ &\geq m_0\epsilon_0 > 2\beta,\end{align*}
a contradiction. So $f$ is continuous from the left on $\R^+$ resp. $\R$. \\
Similarly one shows that $f$ is continuous from the right on $J$, so
$f \in C(J,X)$. (7.1) gives then $\int_t^{t+1} f(s) ds = \int_0^1 f(s)ds + f(t), \ t\in J;$
this shows $f\in C^1(J,X)$ and $f'(t)=f(t+1)-f(t)=f(1)$, or $f(t)=tf(1)$ with $f(0)=0$ (e.g. with \cite[Propositions 1.2.2 and 1.2.3]{lit2}). --\end{proof}
\end{lemma}
\begin{proposition}[De Bruijn] If $I\subset \R$ is an interval and $f:I\to \F \in \{\R,\C\}$ is such that the differences $(\Delta_h f)|(I\cap (I-h))$ are continuous
on $I\cap(I-h)$ for each $h\in\R^+$, then there exist $g\in C(I,\F)$ and $H:\R\to \F$ with (7.1) ($H$ instead of $f$, $J=\R$) with
\begin{equation}
f=g+H \text{ on } I.
\end{equation}
\begin{proof}
This is essentially Theorem 1.3 of \cite[p. 197]{lit17}, though in the proof there some corrections are necessary,
e.g. on p. 204 in Theorem 4.1. additionally translation invariance of the norm is used (on p. 205 below), and on p. 205 above Theorem 4.1 is used with a normed space whose norm is not invariant; this however can be remedied. --  \end{proof}
\end{proposition}
See also \cite[p. 422, Theorem]{lit9}, and the following generalization:
\begin{proposition}
If $J$ is an arbitrary interval $\subset \R$ with non-empty interior, $X$ an arbitrary Banach space, $f:J\to X$ such that
the differences $\Delta_hf|J_h$ are continuous on $J_h$ for all $h\in\R$, $J_h=J\cap (J-h)$, then there exist $G\in C(J,X)$
and $H:\R \to X$ additive (with (7.1) on $\R$) so that
\begin{equation}
f=G+H \text{ on } J.
\end{equation}
\textnormal{The proof is an extension of De Bruijn's proof in \cite{lit17}, with several essential additions, it is too long to reproduce here.}
\end{proposition}
\newpage
\section{Table of properties of some function spaces}
In the following table we collect for easier reference the more important results of this paper and earlier ones,
also a few additional results.

Here always $J=\R$ (though many of the results hold also for $J\ne \R$), $X$ is an arbitrary Banach space (with scalar
field $\F=\C$ if $g$ appears), only $g=e^{it^2}$ is used (again by the results here in many cases more general $\gamma$ are admissible).

(+) or (-) mean some restrictions are needed, e.g. in line 2 for $AP$ the $(P_b)$ and ($L_b$) hold only for $X$ not containing $c_0$ (see \cite[p.120]{lit4}, \cite{lit23}); in line 6 the $(\Gamma)$ holds only for $REC\cap \{f\in UC\mid f(\R)\text{ relatively compact} \}$ by \cite[Theorem 2.6 of Flor]{lit11}.\\

Definitions used which do not appear in the earlier sections:\\

($\Gamma$): $A$ satisfies $(\Gamma)$ means if $f\in A, \ \omega \in \R$ then $\gamma_{\omega} f \in A, \ \gamma_{\omega} (t)=e^{i\omega t}$.\\

($\Delta P$): $A$ satisfies ($\Delta P$) means if $f\in A$, $h\in \R^+$, then $Pf - M_h Pf \in A$ \\

with
	$(Pf)(t):= \int_{\beta}^{t} f(s) ds , \ \beta$ fixed $\in J \ (\beta =0$ if $0\in J$).\\

The proofs can usually be found in the relevant sections §§ 2- 7, for ($\Delta$) see the matrix and the comments at the end of § 6.
Occasionally hints for proofs can be found in the Comments (p. 29) to the line in question.
\newpage
\begin{center}
\textsc{Table of additional results}
\vspace{0.5cm}
\end{center}

\tiny
		\noindent\makebox[\textwidth]{%
    \begin{tabularx}{16.1cm}{|r|r||c|c||c|c|c||c|c||c|c||l|}
    \hline
          &       & Def   & pos.  & unif. & ($\Gamma$) & A$\subset$MA & $(\Delta)$ & ($\Delta P$) & ($P_b$) & ($L_b$) & $\Delta A$  \\
          & A     & p.    & invar. & closed &       &       &       &       &       &       &         \\
          &       &       & (2.6) & (2.8) & ?8   & (2.11) & (2.13) & ?8   & (6.5) & (6.4) & (2.12), (4.19)   \\ \hline
    1     & $C_0$ &  3     & +     & +     & +     & +     & +     & +     & -     & -     & $X_c +C_0+PC_0$   \\
    2     & $AP$  &  3     & +     & +     & +     & +     & +     & +     & (+)   & (+)   & $AP + P(AP)$   \\
    3     & $BAA$ &  7     & +     & +     & +     & +     & -     &       & (+)   & (+)   &         \\ \hline
    4     & $AA$  &  7     & +     & +     & +     & +     & +     & +     &       &       &         \\
    5     & $LAP$ &  7     & +     & +     & (+)   & +     & (+)   & (+)   &       &       &         \\
    6     & $REC$ &  7     & +     & +     & (+)   & +     &       &       &       &       &         \\ \hline
    7     & $PS_+$ &  (3.6)     & +     & +     &       & +     &       &       &       &       &        \\
    8     & $BEM$ &    12   & +     & +     & ?     & +     & -     &       &       &       &         \\
    9     & $BUC$ &    3   & +     & +     & +     & +     & +     & +     & +     & +     & $BUC+PBUC,=UC$   \\ \hline
    10    & $UC$  &    3   & +     & +     & -     & +     & +     & +     & +     & +     & $UC+PUC$   \\
    11    & $BC$  &    3   & +     & +     & +     & +     & +     & +     & +     & +     &         \\
    12    & $L^p$ &    3   & +     & -     & +     & +     & +     & +     & $\text{-}(p\!<\!\infty)$ & $\text{-}(p\!<\!\infty)$ &          \\ \hline
    13    & $\ell^p$ &       & +     & -     & +     &       &      &       &   &  &         \\
    14    & $C_0+X$ &       & +     & +     & -     & +     & +     & +     & -     & -     & $X+Xt+C_0+PC_0$   \\
    15    & $C_0+AP$ &  3     & +     & +     & +     & +     & +     & +     & -     & -     &         \\ \hline
    16    & $C_0+AA$ &  19     & +     & +     & +     & +     & (+)   & (+)   &       &       &         \\
    17    & $EAP$ &   19    & +     & +     & +     & +     & +     & +     & -     & -     &         \\
    18    & $PAP$ &   19    & +     & +     & +     & +     & +     & +     &       &       &         \\ \hline
    19    & $PAA$ &   19    & +     & +     & +     & +     & (+)   & (+)   &       &       &         \\
    20    & $gX$  &   Ex.3.15    & -     & +     & -     & -     & +     & -     & -     & -     & $X_c$   \\
    21    & $gAP$ &   Ex.3.14    & +     & +     & +     & -     & -     & -     & -     & -     & $X_c+gAP$   \\ \hline
    22    & $gAA$ &   Ex.3.14    & +     & +     & +     & (-)   &  (-)      & -     & -     & -     &         \\
    23    & $gBUC$ &   (4.15)    & +     & +     & +     & +     & +     & +     & -     & -     & $X_c+gBUC+PC_0$   \\
    24    & $gUC$ &   (4.15)    & -     & +     & -     & +     & +     & +     & -     & -     & $\Delta(gBUC)$  \\ \hline
    25    &       &       &       &       &       &       &       &       &       &       &         \\
    26    & $X+gX$ &  (6.42)     & -     & +     & -     & -     & +     & -     & -     & +     & $X_c+Xt$   \\
    27    & $X+gAP$ &       & +     & +     & -     & -     & -     &       & -     & +     & $X_c+Xt+gAP$   \\ \hline
    28    & $X+gAAP$ &       & +     & +     & -     & +     & +     & +     &       &       &$X_c \!+ \! Xt  \!+ \! gAAP \! + \! PC_0   $         \\
    29    & $X+gBUC$ &       & +     & +     & -     & +     & +     & +     & -     & -     & $X\!+\! Xt\! +\! gBUC\! +\! PC_0$   \\
    30    & $X+gUC$ &       & -     & +     & -     & +     & +     & +     & -     & -     & $\Delta (X +gBUC)$   \\ \hline
    31    & $C_0+gX$ &   (6.42)    & -     & +     & -     & +     & +     & +     & -     & -     & $\Delta C_0$   \\
    32    & $C_0+gAP$ &   29,C.32    & +     & +     & +     & +     & +     & +     & -     & -     & $X + gAAP + PC_0$        \\
    33    & $AP+gX$ &  (6.42)     & -     & +     & -     & -     & +     & -     & -     & (+)   & $\Delta AP$   \\ \hline
    34    & $AP+gAP$ &       & +     & +     & +     & -     & -     & -     & -     & -     &         \\
    35    & $AP+gAAP$ &       & +     & +     & +     & +     & +     & +     & +     & +     &         \\
    36    & $AP+gBUC$ &       & +     & +     & +     & +     & +     & +     & -     & -     &        \\ \hline
    37    & $AP+gUC$ &   25    & -     & +     & -     & +     & +     & +     & +     & +     &         \\
    38    & $AAP+gX$ &   (6.42)    & -     & +     & -     & +     & +     & +     & -     & -     &         \\
    39    & $BUC+gX$ &   (6.42)    & -     & +     & -     & +     & +     & +     & +     & +     & $BUC+PBUC$   \\ \hline
    40    & $BUC+gAP$ &       & +     & +     & +     & +     & +     & +     & +     & +     &         \\
    41    & $BUC+gBUC$ &       & +     & +     & +     & +     & +     & +     & +     & +     & $BUC+gBUC+PBUC$   \\
    42    & $BUC+gUC$ &   25    & -     & +     & -     & +     & +     & +     & +     & +     &  $\Delta (BUC + gBUC)  $      \\ \hline
    43    & $UC+gX$ &   (6.42)    & -     & +     & -     & +     & +     & +     & +     & +     & $UC+PUC$   \\
    44    & $UC+gAP$ &       & +     & +     & -     & +     & ?     & +     & +     & +     &         \\
    45    & $UC+gBUC$ &       & +     & +     & -     & +     & ?     & +     & +     & +     &         \\ \hline
    46    & $UC+gUC$ &   25    & -     & +     & -     & +     & ?     & +     & +     & +     &         \\
    47    & $C_0 + \gamma_{\omega}X$ & (2.2)      & +     & +     & -     & +     & +     & +     & -     & -     & $C_0+\gamma_{\omega}X+X+PC_0$   \\
    48    &       &       &       &       &       &       &       &       &       &       &         \\ \hline
    49    &       &       &       &       &       &       &       &       &       &       &         \\
    50    &       &       &       &       &       &       &       &       &       &       &         \\
    \hline
    \end{tabularx}%
		}
 \label{tabelle}
\normalsize

For ($\Delta$) see also the matrix on p. 25.\\
\newpage
\begin{center}
\textbf{Comments}\\
(The numbers refer to the corresponding row of the table above. )
\end{center}
1. $C_0 = gC_0 ; (\Delta)$: Prop. 4.1; $\Delta C_0$ : (4.19).\\
4. [9, Prop. 3.5 (ii)] .\\
5. [9, Prop. 3.2, 3.5(i)].\\
6. Not linear.\\
7. Not linear.\\
8. Example 4.8, (4.22).\\
9. Prop. 6.4 (b), after (6.4); Ex. 4.7.\\
10. $ UC = BUV+PBUC : UC \subset \Delta BUC$, see 9.\\
11. [8, Ex. 3.5], Prop. 6.4(b).\\
12. Ex. 3.4, [8, Prop. 3.4].\\
15. $C_0+AP = : AAP$, p. 3.\\
18. $PAP = \text{Pseudo } AP$, \cite{lit29}.\\
20. Ex. 3.15, Ex. 4.4.\\
21. Ex. 3.14. (3.9).\\
23. (4.15).\\
24. (4.15), Prop. 6.14 .\\
26. Matrix p. 24, p. 25.\\
27. - 30. p. 24/25.,  Prop. 6.14.\\
31. Ex.s 6.11.\\
32. $C_0 + gAP = gAAP$: See Proof of Ex. 6.2.\\
33. Ex.s 6.11, Prop. 6.10.\\
34. p. 24/25.\\
35. $AP+gAAP=AAP+gAP=AAP+gAAP$, p. 24/25.\\
36. $AP+gBUC=AAP+gBUC$, p. 24/25.\\
37. - 42. P. 24/25.\\
38., 39., 43. Ex.s 6.11, Prop. 6.10.\\
40. $BUC+gAP=BUC+gAAP$, p. 24/25, Prop. 6.4(b).\\
41. Theorem 6.5, Cor. 6.6, Ex.s 6.8, Cor. 6.9.\\
42. p.24/25, Prop. 6.14.\\
43. Ex.s 6.11, Prop. 6.10.\\
44., 45., 46.  $(L_b)$ : Porada inequ.\\
47. Prop. 4.1, (4.19).\\

\newpage
\section{Open questions} 
\begin{compactenum}[1.]
\item Does there for $J=\R_+$ exist a non-trivial $\lambda$-class satisfying the Loomis condition ($L_b$) but $\ne BUC(J,X)$?\\
($\lambda$-class: linear closed subspace of $BUC$ with ($\Gamma$) (see § 8) and in a suitable sense invariant \cite[p. 117]{lit4}.)
\item Does there exist a simpler proof for ($\Delta$) for $BUC + \gamma BUC$ than ours in §6?
\item Does there exist a simple example of an $A$ with $A\subset MA $ but without ($\Delta$), and if possible linear, invariant, uniformly closed, $\subset BC$?\\
(See example 4.8, and Q4.)
\item Is there an extension of Proposition 4.1 which explains why reasonable $A\subset MA$ and $\subset BC$ but without ($\Delta$) are rare?
\item Is there a simple proof of ($\Delta$) for $gAAP$? Extension to $g_{\omega, r}AAP, \gamma AAP$?
\\ (True by an unpublished result of B. Basit and the author 2002.)
\item Is there a simple example of a (linear, invariant, uniformly closed) $A\ne\{0\}$ with $A\subset MA$ and ($P_b$), but without ($L_b$)?\\
(Definitions: (6.4), (6.5); the converse situation holds e.g. for $X_c+gX, $ or $AP+gX$; see Q7.)
\item Is there a theorem which explains why reasonable $A$ with ($P_b$) but not ($L_b$) are rare? (see Q6., Q4.)
\item Does $UC+gUC$ satisfy ($\Delta$)?
\item Which of the other two ``?'' in the $\Delta$-matrix in § 6 can be replaced by a ``+''? ($ UC+gAP=UC+gAAP$, Q8.)
\item Examples of $U, V\subset X^J$ with $MU+MV$ strictly $\subset M(U+V)$?\\
(See \cite[Prop.5.1 (ii) p. 38]{lit6}?)
\item Can one improve Example 3.18 to $g_{\omega,r}UC\subset M(g_{\omega,r}U)$ with $U$ smaller than $BUC$?
\item Can one extend Examples 3.18/3.19 to $\gamma\ne g_{\omega,r}$, but still with $UC\cap (\gamma UC)\subset C_0, \gamma UC \not\subset C_0$?
\item Can one improve (3.8) to ($\sum_0^nM^m(\gamma A))\cap (\sum_{n+1}^{\infty}M^m(\gamma A))=\{f=0 \text{ a.e.}\}, n\in \N$?
\item Inclusion relations between $\cup_0^{\infty} M^m (gAP)$ and $\mathscr{D}'_{gAP}\cap\Lloc$?\\
(See Prop. 3.7 (i),(ii).)
\item Example of $A$ with ($\Delta P$) but not ($\Delta$) (and $A$ linear, positive invariant, uniformly closed, $A\subset MA$)? (See (4.20); $A$ with ($\Delta$) but not ($\Delta P$): $gX, X_c + gX; (\Delta P)$: § 8.)
\item Can one extend ``$gAAP$ has $(\Delta)$'' to ``$\gamma(C_0+U)$ has ($\Delta$)'' with $U\subset PS_+ $ (and ($\subset (B)UC)$ with ($\Delta$)? See Q5.
\item Can one weaken $\inf_J |\gamma| >0$ resp. $\liminf_{|t|\to\infty}|\gamma|>0$ to $\liminf_{|t|\to\infty}|\gamma|>0$ resp. $\limsup_{|t|\to\infty}|\gamma|>0$ or $\gamma(t)\ne 0$ (for $t\in J$ or $|t|\geq n_0$), e.g. in Corollaries 5.2/6.6, (5.4), (6.7)?
\item Construct a complete Venn diagramm for the $\gamma$-properties $O_1,O_2, \gamma \in MC_0$, Porada (5.9) holds, $UC\cap \gamma UC \subset C_0$,
with/without $\norm{\gamma}_{\infty} < \infty $ and/or $\liminf_{|t|\to\infty}|\gamma|$\\$>0$ and/or $\gamma \in C(J,\F)$.
\item Example of a $\gamma \ne g_{\omega, r}$ which satisfies all the assumptions in (6.7) or (6.41) or for which $\gamma UC\subset M(\gamma UC)$?
\item Example of a non-trivial $\gamma\ne g_{\omega,r}$ for which $BUC +\gamma BUC$ satisfies ($\Delta$) and is uniformly closed?
\item Can one reduce/simplify the assumptions in (6.7) or (6.41), or are all the assumptions there necessary for Theorem 6.5 resp. Corollary 6.6?
\item Is Theorem 6.5 true, if (6.7) is assumed only for the $J$ in Theorem 6.5? (See the proof of Corollary 6.6, and Q23.)
\item If (6.7) is asumed only for some $J\ne \R$, can one extend $\gamma $ to $\R$ so that (6.7) holds for this extension and $J=\R$?
\item Do there exist real-valued $\gamma $ satisfying (6.7)? (Not possible for (6.41).)
\item Are there $r>2$ so that $BUC + g_{1,r} BUC $ satisfies $(L_b)$ or ($\Delta$)?
\item Does $BUC(J,X) + \cos(t^2) BUC(J,X)$ satisfy ($\Delta$) or ($L_b$) ($\F=\R $ or $\C$)?
\item Example of $\gamma$ with $(B)UC \cap \gamma (B) UC)\subset C_0, $ but $\gamma \notin C_0$, $\gamma $ not satisfying $O_1$?
\item Do $A\subset MA $ and $(\Delta)$ hold for $A=g_{1,r} UC(J,X), 0<r<1$?
\item If $UC+gUC$ does not satisfy ($\Delta$), does $A+gA$ satisfy ($\Delta$) for $A=AP(J,\C) \cdot UC(J,X)$, or more gernerally $A+\gamma A$ for $A=V\cdot UC$, with some $V (\subset BUC)$ (see Q8.)?
\item Is $g=e^{it^2}$ Maak-ergodic (see (4.22))?
\end{compactenum}


\noindent
Mathematisches Seminar\\
Universität Kiel\\
D 24098 Kiel\\
Germany\\
hans@guenzler.de
\end{document}